[1994/12/01]
\documentclass[12pt]{article}
\title{Crystallography and Riemann Surfaces}
\author{Veit Elser\\
Department of Physics\\
Cornell University}

\input epsf         
\epsfverbosetrue    

\usepackage{amsmath,amsthm}

\usepackage{graphics}

\usepackage{amsfonts}

\chardef\bslash=`\\ 

\newtheorem{thm}{Theorem}[section]
\newtheorem{cor}[thm]{Corollary}
\newtheorem{lem}[thm]{Lemma}

\theoremstyle{definition}
\newtheorem{defn}{Definition}[section]

\theoremstyle{remark}

\newtheorem*{notation}{Notation}

\newcommand{\A}{\mathcal{A}}

\newcommand{\Poly}{\mathcal{P}} 
 
 \newcommand{\F}{\mathcal{F}} 
 \newcommand{\Surf}{\mathcal{S}} 
\newcommand{\T}{\mathcal{T}}

\DeclareMathOperator{\order}{O}

\DeclareMathOperator{\Ker}{Ker}
\DeclareMathOperator{\LCM}{LCM}
\DeclareMathOperator{\Aut}{Aut}
\DeclareMathOperator{\GCD}{GCD}
\DeclareMathOperator{\Rk}{rk}
\DeclareMathOperator{\GL}{GL}
\DeclareMathOperator{\SL}{SL}

\newcommand{\eval}[2][\right]{\relax
  \ifx#1\right\relax \left.\fi#2#1\rvert}

\begin{document}
\maketitle

\begin{abstract}
The level set of an elliptic function is a doubly periodic point 
set in $\mathbb{C}$.  To obtain a wider spectrum of point sets, we 
consider, more generally, a Riemann surface $\Surf$ immersed in 
$\mathbb{C}^{2}$ and its sections (``cuts'') by $\mathbb{C}$.  We give 
$\Surf$ a crystallographic isometry in $\mathbb{C}^{2}$ by defining a 
fundamental surface element as a conformal map of triangular domains 
and $\Surf$ as its extension by reflections in the triangle edges. 
Our main result concerns the 
special case of maps of right triangles, with the right angle being a 
regular point of the map.  For this class of maps we show that only 
seven Riemann surfaces, when cut, form point sets that are discrete in 
$\mathbb{C}$.  Their isometry groups all have a rank-four lattice 
subgroup, but only three of the corresponding point sets are doubly 
periodic in $\mathbb{C}$.  The remaining surfaces form quasiperiodic 
point sets closely related to the vertex sets of quasiperiodic 
tilings.  In fact, vertex sets of familiar tilings are recovered in 
all cases by applying the construction to a piecewise flat 
approximation of the corresponding Riemann surface.  The geometry of 
point sets formed by cuts of Riemann surfaces is no less ``rigid'' than the 
geometry determined by a tiling, and has the distinct advantage in 
having a regular behavior with respect to the complex parameter which 
specifies the cut.
\end{abstract}

\renewcommand{\sectionmark}[1]{}

\section{Introduction}

Crystallography is concerned with point sets in $\mathbb{R}^{n}$ that are 
discrete and distributed more or less uniformly. In ``classical'' 
crystallography, periodicity was imposed as well, but this 
restriction is not considered as fundamental in the ``modern'' era. With 
the discovery of intermetallic quasicrystals \cite{AlMn} in the early 
1980s, it became clear that there exist aperiodic point sets that 
share a basic property with periodic point sets.  This property should 
really be associated with a distribution: in this case, the 
distribution formed by placing a Dirac delta at each point of the set.  
In these terms, the class of aperiodic sets singled out by 
crystallography is characterized by the property that the Fourier 
transform of the corresponding distribution has support on a lattice 
\cite{Mermin}.  The rank of this ``Fourier lattice'' equals the dimension 
of space for periodic sets, and exceeds it (but is still finite), in 
the case of \emph{quasiperiodic} sets.  The vertex set of the Penrose 
tiling of the plane \cite{Penrose} is a familiar example of a 
quasiperiodic set.

The standard construction of quasiperiodic sets $\A$ begins by embedding 
$\mathbb{R}^{n}=Y$ in a larger Euclidean space, 
$\mathbb{R}^{m+n}=X\times Y$. Into $X\times Y$ one then immerses a 
smooth $m$-manifold $\Surf$ that is (\textit{i}) transversal to $Y$, 
and (\textit{ii}) invariant under the action of a lattice $\Lambda$ 
generated by $n+m$ linearly independent translations in $X\times Y$.  
Point sets $\A\subset Y$ are obtained as sections (``cuts'') of 
$\Surf$ by spaces parallel to $Y$.  More formally, in terms of the 
standard projections
\begin{equation}
\begin{align}
\pi_{X}&\colon \Surf\to X\\
\pi_{Y}&\colon \Surf\to Y,
\end{align}
\end{equation}
the section of $\Surf$ at $x\in X$ is the set
\begin{equation}\label{pointSets}
\A(x)=\pi_{Y}\circ{\pi_{X}}^{-1}(x).
\end{equation}
Periodicity or quasiperiodicity of $\A(x)$ is determined by the rank 
of the lattice $\Lambda_{Y}=\Lambda\cap Y$. Since the generators of 
$\Lambda$ were assumed to be linearly independent, 
$\Rk{(\Lambda_{Y})}\leq n$. Quasiperiodicity corresponds to 
$\Rk{(\Lambda_{Y})}<n$, with complete absence of periodicity 
characterized by $\Lambda_{Y}=\{0\}$. An important motivation for 
constructing quasiperiodic sets, in this context, is the fact that symmetry 
groups that cannot be realized by periodic point sets in 
$\mathbb{R}^{n}$, \emph{can} 
be realized by periodic surfaces in $\mathbb{R}^{m+n}$.

Transversality and periodicity are relatively mild restrictions on 
the manifold $\Surf$, called the ``atomic surface'' by physicists. A 
further restriction, one which leads to point sets called ``model 
sets'' \cite{Moody}, is to require that $\Surf$ is the $\Lambda$-orbit of 
a polytope in $X$.  The algorithm which constructs $\A(x)$ from such 
$\Surf$ naturally leads to the terminology ``window'' or ``acceptance 
domain'' for the corresponding polytopes.  Model sets can always be 
organized into finitely many tile shapes, and, because of this 
simplicity, have dominated the study of quasiperiodic sets.

A different viewpoint on the construction of $\Surf$, pioneered by 
Kalugin \cite{Kalugin} and Katz \cite{Katz}, emphasizes the continuity properties of $\A(x)$ 
with respect to $x$. Consider in more detail the construction of a 
model set: $\Surf=\Poly+\Lambda$, where $\Poly\subset X$ is a 
polytope. Now, if $x\in \Poly+\pi_{X}(\lambda)$ for some 
$\lambda\in\Lambda$, then $y=\pi_{Y}(\lambda)\in \A(x)$. But now 
consider what happens when $x$ crosses the boundary of 
$\Poly+\pi_{X}(\lambda)$. As $x$ ``falls off the edge of the earth'', 
the corresponding point $y$ in the point set $\A(x)$ disappears. By 
the same process, of course, points can spontaneously appear ``out of 
thin air''.  To gain control over these processes, Kalugin \cite{Kalugin} 
and Katz \cite{Katz} advocated a restriction on $\Poly$, in relation to 
$\Lambda$, such that whenever $x$ falls off the edge of one polytope, 
$\Poly+\pi_{X}(\lambda)$, it falls within another, say 
$\Poly+\pi_{X}(\lambda')$.  This restriction corresponds mathematically 
to the statement that the boundaries of the disconnected components of 
$\Surf=\Poly+\Lambda$ can be ``glued'' together to form a topological 
manifold without boundary.

In the process of restoring transversality to the glued complex of 
polytopes one encounters the problems addressed by singularity theory. 
The map $\pi_{X}$ should now be a smooth (but not necessarily 1-to-1) 
map of $m$-mainfolds. 
In the trivial 
situation, when $\pi_{X}$ has no singularities, $\Surf$ must be 
diffeomorphic to a collection of hyperplanes. This is the situation explored by 
Levitov \cite{Levitov} for point sets in two and three dimensions and various 
symmetry groups. 

When $\Surf$ is a generic 2-manifold, we have the classic result of 
Whitney \cite{Whitney} that the stable singularities of smooth maps, such as 
$\pi_{X}$, are folds and cusps, having respectively codimension one 
and two.  Because the cusp is always accompanied by two folds, the 
locus of singular values of $\pi_{X}$ consists of curves.  The space 
$X$ is thus populated by singular curves such that whenever $x$ 
crosses a curve, a pair of points in $\A(x)$ merge and annihilate.  
One motivation for the present work was the desire to eliminate this 
point-merging singularity to the greatest extent possible.

By giving the 2-manifold $\Surf$ a complex structure, and identifying 
$X$ with the complex plane, we impose additional regularity by 
insisting that $\pi_{X}$ is locally holomorphic. The singularities of $\pi_{X}$ 
will then be isolated points. A construction that naturally
leads to a $\pi_{X}$ with this property is to let $\Surf$ be (locally) 
the graph of a holomorphic function, $f\colon X\to Y$. Globally this 
corresponds to a Riemann surface $\Surf$ immersed in $\mathbb{C}^{2}=X\times 
Y$ and having an atlas of compatible charts in $X$. The other 
ingredient needed by our construction is some way to guarantee 
that $\Surf$ is invariant with respect to a lattice $\Lambda$. We meet 
this challenge by using conformal maps between triangles to define a 
fundamental graph of $\Surf$. Schwarz reflections in the triangle 
edges extend this graph and generate the isometry group of $\Surf$. 
For appropriate choices of triangles, the isometry group has a lattice 
subgroup with the desired properties.

In the second half of this paper we classify a subset of all 
Riemann surfaces generated by conformal maps of triangles. This 
subset is characterized by the property that the conformal map is regular at one 
vertex of the triangles and that the edges at this vertex make the 
largest possible angle, $\pi/2$. With the only other restriction being 
that the corresponding point sets $\A(x)$ are 
discrete in $Y$, one arrives at a set of seven surfaces. Four of these are 
quasiperiodic. The point set 
obtained from a section of one of them is shown in Figure 1. Also 
shown in Figure 1 is a much studied model set \cite{pent}: a tiling of 
boats, stars, and jester's-caps (whose vertices coincide with a subset 
of the Penrose-tiling vertex set).  The point set determined by the 
Riemann surface can be said to be approximated by the vertex set of 
the tiling by a systematic process that renders the Riemann surface 
piecewise flat.  The other surfaces obtained in our partial 
classification, when flattened, also produce familiar tilings (Fig.  
3).

\begin{figure}
\centerline{\epsfbox{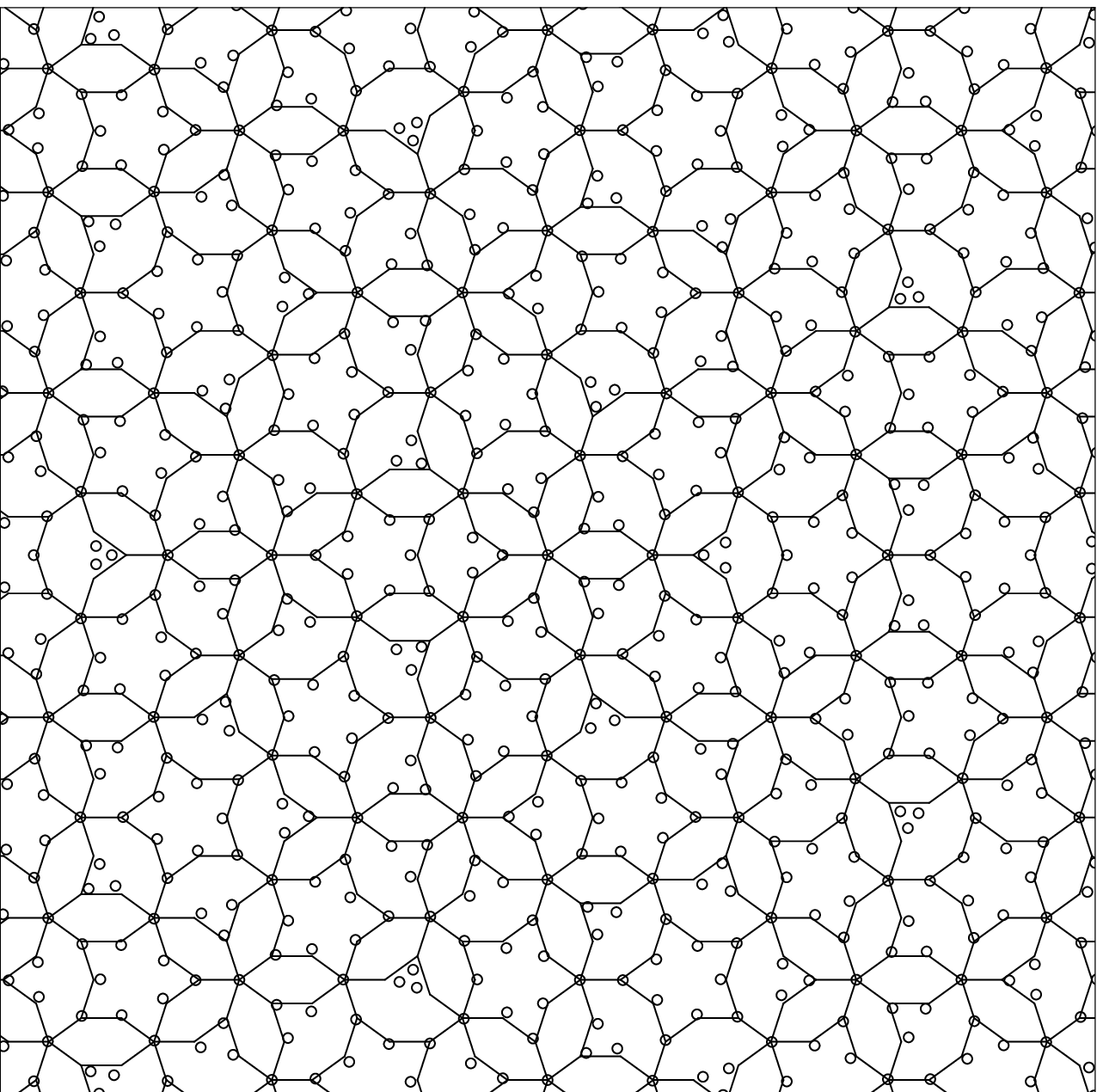}}
\caption{Point set (small circles) given by the section of a Riemann 
surface compared with a tiling of boats, stars, and jester's caps \cite{pent}.}
\end{figure}

\newpage

\section{Riemann surfaces generated by conformal maps}

\subsection{Immersed Riemann surfaces}

We consider Riemann surfaces as analytically continued holomorphic 
functions, interpreted geometrically as surfaces immersed in 
$\mathbb{C}^{2}$.  Our treatment follows closely the notation and 
terminology of Ahlfors \cite{Ahlfors}.

\begin{defn}
A \emph{function element} $F=(U,f)$ consists of a domain $U\subset 
\mathbb{C}$ and a holomorphic function $f\colon U\to \mathbb{C}$.
\end{defn}

\begin{defn}
Function elements $F_{1}=(U_{1},f_{1})$ and $F_{2}=(U_{2},f_{2})$ are \emph{direct 
analytic continuations} of each other iff $V=U_{1}\cap U_{2}\neq\emptyset$
and $f_{1}=f_{2}$ when restricted to $V$.
\end{defn}

\begin{defn}
The \emph{complete, global analytic function determined by function element}
$F_{0}=(U_{0},f_{0})$ is the maximal collection of function elements 
$\F$
such that for any $F_{i}\in \F$ there exists a chain of 
function elements $F_{0},\ldots,F_{i}$, all in $\F$, 
with every link in the chain a direct analytic continuation.
\end{defn}

Up to this point the set of function elements comprising a complete, 
global analytic function $\F$ only possesses the discrete topology, 
where $F_{1}\cap F_{2}=\emptyset$ whenever $F_{1}$ and $F_{2}$ are 
distinct elements of $\F$. By refining this topology we can identify 
$\F$ with a surface and, ultimately, a Riemann surface.
Consider a pair of function elements in $\F$, $F_{1}=(U_{1},f_{1})$ and 
$F_{2}=(U_{2},f_{2})$. In the refined topology we define the 
intersection by
\begin{equation}
F_{1}\cap F_{2}=\left\{
\begin{array}{ll}
F_{3}=(U_{3},f_{3}) & \parbox[t]{1.75in}{if $F_{1}$ and $F_{2}$ are related by 
direct analytic continuation,} \\
& \\
\emptyset & \parbox[t]{1.75in}{otherwise,}
\end{array}\right.
\end{equation}
where $U_{3}=U_{1}\cap U_{2}$ and $f_{3}$ is $f_{1}=f_{2}$ restricted 
to $U_{3}$. It is straightforward to check that this defines a valid 
topology. Moreover, the projection $\pi\colon \F\to \mathbb{C}$ given 
by
\begin{equation}
\pi\colon (U,f)\mapsto U,
\end{equation}
provides the complex charts that identify $\F$ with a Riemann surface.

Throughout the rest of this paper we will mostly be interested in 
Riemann surfaces immersed in $\mathbb{C}^{2}$.
\begin{defn}
Let $\F$ be a complete, global analytic function. The \emph{immersed 
Riemann surface $\Surf$ corresponding to $\F$} is the image of the 
immersion $\Psi\colon \F\to \mathbb{C}^{2}$ given by
\begin{equation}
\Psi\colon (U,f)\mapsto \{ (x,f(x))\colon x\in U \}.
\end{equation}
\end{defn}

\begin{notation} We denote the first component of $\mathbb{C}^{2}$ by 
$X$, the second by $Y$.
\end{notation}

If 
we restrict the immersion $\Psi$ to a single function 
element, $F_{0}=(U_{0},f_{0})$, we obtain the \emph{graph}
\begin{equation}
\Surf_{0}=\{ (x,f_{0}(x))\in X\times Y\colon x\in U_{0}\}\label{graph}
\end{equation}
Thus $\Surf_{0}$ represents a piece of $\Surf$ and in fact determines all of 
$\Surf$; $\Surf$ is 
connected because every pair of function elements in a complete 
global analytic function is related by a chain of direct analytic 
continuations.  $\Surf$ is the \emph{completion} of $\Surf_{0}$.
\begin{notation}
We write $[\Surf_{0}]$ to denote the completion of the graph $\Surf_{0}$.
\end{notation}
Since all subsequent references to ``Riemann surface'' will be as 
a surface immersed in $X\times Y$, we drop the qualifier ``immersed'' below.
We also omit the term ``complete'', since the only instances of 
incomplete surfaces, graphs, will always be identified as such. Given a 
Riemann surface $\Surf$, we will frequently make use of the 
projections
\begin{equation}
\begin{align}
\pi_{X}&\colon \Surf\to X,\\
\pi_{Y}&\colon \Surf\to Y.
\end{align}
\end{equation}

The historical construction of Riemann surfaces we have followed can 
be criticized for its inequivalent treatment of the spaces $X$ and $Y$. 
We can correct this fault by insisting that the functions $f$ 
appearing in the function elements $(U,f)$ are not just holomorphic 
in their respective domains $U$, but \emph{conformal} 
(holomorphic with holomorphic inverse). 
The graph \eqref{graph} could then be equally written as
\begin{equation}
\Surf_{0}=\{ ({f_{0}}^{-1}(y),y)\in X\times Y\colon y\in V_{0}\},
\end{equation}
where $V_{0}=f_{0}(U_{0})$. If this ``inversion'', or interchange of 
$X$ with $Y$, is to work for all function elements $(U,f)$, 
then one must remove all points $x_{0}\in U$, where $f$ behaves 
locally as $f(x)-f(x_{0})=c(x-x_{0})^{m}+\dotsb$, with $m>1$. These correspond to branch 
points of the map $\pi_{Y}$. Conversely, had we begun with the 
inverted function elements our Riemann surface would have \emph{included} 
branch points of $\pi_{X}$, \textit{i.e.} points where $f$ is singular.
In keeping with tradition we augment our definition of a Riemann 
surface $\Surf$
to \emph{include} all points $(x_{0},y_{0})$ where $\Surf$ behaves 
locally like the algebraic curve $(y-y_{0})^{n}=c(x-x_{0})^{m}$, 
where $m$ and $n$ are positive integers.
\begin{defn}
A point $(x,y)\in \Surf$ is \emph{regular} if the corresponding complete, global 
analytic function
contains a function 
element $(U,f)$, with $x\in U$ and $f$ conformal at $x$. A point which is 
not regular is \emph{singular}.
\end{defn}

\subsection{Transformations}

Two transformations of Riemann surfaces will be needed in our 
discussion of symmetry properties. These are defined in terms of 
their action on the spaces $X$ and $Y$ and induce a transformation on 
Riemann surfaces as subsets of $X\times Y$. Let $(x,y)$ be a 
general point in $X\times Y$ and define  
the following transformations:
\begin{align}
\tau(a,b;c,d)&: (x,y)\mapsto (a x+b ,c y+d )\\
\sigma &: (x,y)\mapsto (\bar{x},\bar{y})
\end{align}
Transformation $\tau$ (for complex constants $a$, $b$, $c$, and $d$) 
is the general bilinear map while $\sigma$ corresponds 
to Schwarz reflection (componentwise complex conjugation).
\begin{lem}
If $\Surf$ is a Riemann surface and $T$ is either of the transformations 
$\tau$ or $\sigma$, then $T\Surf$ is again a Riemann surface.
\end{lem}
\begin{proof}
Write $T(x,y)=(T_{X}x, T_{Y}y)$ where $T_{X}$ and $T_{Y}$ are just 
maps of the complex plane. Since $\Surf$ corresponds to a 
complete global analytic function 
$\F$, we need to verify 
that $T\Surf$ corresponds to some other complete global analytic 
function $\F_{T}$. From our definitions we see that $\F_{T}$ is 
obtained from $\F$ by substituting each function element $F=(U,f)\in \F$
by $F_{T}=(T_{X}U,T_{Y}\circ f\circ T_{X}^{-1})$. It is 
easily checked that $T_{X}$ is open and $T_{Y}\circ f\circ 
T_{X}^{-1}$ is holomorphic for both of the transformations 
being considered. Thus $F_{T}$ remains a valid function element. One 
also verifies that the direct analytic continuation relationships 
among function elements are unchanged by these transformations.
\end{proof}

\begin{cor}\label{moveTinside}
If $\Surf_{0}$ is a graph and $T$ is either of the transformations 
$\tau$ or $\sigma$, then $[T\Surf_{0}]=T[\Surf_{0}]$.
\end{cor}

Two special transformations are rotations and translations, for 
which we introduce the following notation:
\begin{align}
r(\theta,\phi)&=\tau(e^{i\theta},0;e^{i\phi},0)\\
t(u,v)&=\tau(0,u;0,v).
\end{align}

More generally, transformations $T\colon X\times Y\to X\times Y$ 
which act isometrically on the spaces $X$ and $Y$ are 
just the products of Euclidean motions in $X$ and $Y$. 
Isometries of Euclidean spaces normally include reflections; to 
preserve the structure of the immersed Riemann surface, however, any
reflection in $X$ (complex conjugation) must be accompanied by a 
reflection in $Y$.

\begin{defn}
The \emph{isometries of} $X\times Y$ is the group of transformations 
generated by $\sigma$, $r(\theta,\phi)$, and $t(u,v)$.
\end{defn}

In what follows we use the term ``isometry'' only in this sense.
Isometries which preserve a Riemann surface $\Surf$ are called 
isometries of $\Surf$ and form a group. The maximal group of isometries 
is called the isometry group of $\Surf$.

\begin{defn}
The group of \emph{proper isometries} of $X\times Y$ is the normal subgroup of 
isometries generated by $r(\theta,\phi)$ and $t(u,v)$. Any element of 
the coset, $\sigma g$, where $g$ is a proper isometry, is called a 
\emph{Schwarz reflection}. 
\end{defn}

\subsection{Surfaces generated by conformal maps of triangles}\label{conformalMapTriangle}

We now focus on the class of Riemann surfaces determined 
by graphs which solve a purely geometrical problem: 
the conformal map between two bounded triangular regions, $P\subset X$ and 
$Q\subset Y$.
The Riemann mapping theorem \cite{Ahlfors} asserts there is a three-parameter 
family of conformal maps $f\colon P\to Q$ that extend to 
homeomorphisms of the closures $\Bar{P}$ and $\Bar{Q}$.  To fix these 
parameters we require that the three vertices of $\Bar{P}$ map to 
the vertices of $\Bar{Q}$.  This defines the graph
\begin{equation}
P|Q=\{(x,f(x))\colon x\in P\},
\end{equation}
and a corresponding Riemann surface $[P|Q]$. The closure of $P|Q$ 
is defined analogously and is written $\Bar{P}|\Bar{Q}$.  One of the 
main benefits of using a conformal map of triangles to determine a 
Riemann surface $\Surf$ is that its isometry group can be understood 
simply in terms of its action on a partition of $\Surf$ into 
\emph{tiles}.

Just as $\Bar{P}$ can be decomposed into an interior $P$, edges which 
bound $P$, and vertices which bound each edge, there is a 
corresponding cell decomposition of the graph $\Bar{P}|\Bar{Q}$.  For 
example, if $P_{1}$ is one vertex of $\Bar{P}$, and $f(P_{1})=Q_{1}$ is 
its image in $\Bar{Q}$, then we use the symbol $P_{1}|Q_{1}$ to 
represent the corresponding vertex of $\Bar{P}|\Bar{Q}$.  Each vertex of 
$\Bar{P}|\Bar{Q}$ is associated with two angles, a vertex angle of $P$ 
and the corresponding vertex angle of $Q$.  Let the three angle pairs 
be $\alpha_{i},\beta_{i}$, $i=1,2,3$.  If $\alpha_{i}=\beta_{i}$ for 
all $i$, then $P$ is similar to $Q$ and $f$ is just a linear map.  
Because the corresponding Riemann surface would be trivial (a plane) 
we exclude this case.  It is impossible to have 
$\alpha_{i}\neq\beta_{i}$ for just one $i$ since then the angle sum 
could not be $\pi$ in both triangles.  Thus we must have at least two vertices with 
unequal angles.  At these vertices $f$ fails to be conformal.  Any 
vertex of $\Bar{P}|\Bar{Q}$ where the corresponding angles in $P$ and 
$Q$ are unequal will be called a \emph{singular vertex}.  The singular 
vertices of $\Bar{P}|\Bar{Q}$ are the only singular points of 
$\Bar{P}|\Bar{Q}$.

The edges of $\Bar{P}|\Bar{Q}$ (associated with each pair of vertices 
$ij=12,13,23$) are effectively the generators of the isometry group of 
$[P|Q]$.  By $P_{12}|Q_{12}$ we mean the graph given by the 
restriction of $f$ to the edge $P_{12}$ of $\Bar{P}$ with image an edge 
$Q_{12}$ of $\Bar{Q}$.  Consider the triangles $P'$ and $Q'$ obtained 
from $P$ and $Q$ by reflection in these edges.  The graph $P'|Q'$, 
determined by the conformal map $g\colon P'\to Q'$, is clearly related 
to $P|Q$ by an isometry of $X\times Y$, a Schwarz reflection which we 
call $\sigma_{12}$.  Because $\sigma_{12}$ fixes every point of 
$P_{12}|Q_{12}$, we have that $f(x)=g(x)$ for all $x\in P_{12}$.  A 
basic result from complex analysis then tells us that the function 
elements $(P,f)$ and $(P',g)$ are related by analytic continuation.  
Thus $[P|Q]= 
[P'|Q']=[\sigma_{12}(P|Q)]=\sigma_{12}[P|Q]$, by Corollary \ref{moveTinside}.
\begin{defn}
The group $G$ generated by the 
Schwarz reflections $\sigma_{ij}$ which fix the three edges of 
a triangular graph $P|Q$ is called the \emph{edge group} of $P|Q$. The 
edge group of $P|Q$ is a subgroup of the isometry group of $[P|Q]$.
\end{defn}
In order 
to show that the edge group of a triangular graph is the \emph{maximal} 
isometry group, we first need to 
refine the sets on which these groups act.

\begin{notation}
The symbol $\Check{\Surf}$ corresponds to the Riemann surface $\Surf$ 
whose singular points have been removed.
\end{notation}
\begin{defn}
A \emph{real curve} of the Riemann surface $\Surf$ is any curve $\Gamma\subset 
\Check{\Surf}$, homeomorphic to 
$\mathbb{R}$, and pointwise invariant with respect to a Schwarz 
reflection.
\end{defn}
Since both $\pi_{X}\colon\Check{\Surf}\to X$ and 
$\pi_{Y}\colon\Check{\Surf}\to 
Y$ are immersions, the map $\pi_{Y}\circ{\pi_{X}}^{-1}\colon 
\pi_{X}(\Gamma)\to \pi_{Y}(\Gamma)$ is an immersion as well. Thus it 
makes sense to use our graph notation, $\Gamma=\gamma|\delta$, for real 
curves, where $\gamma=\pi_{X}(\Gamma)$ and $\delta=\pi_{Y}(\Gamma)$.  
A real curve $\gamma|\delta$ is geometrically no different from the 
edge of a triangular graph; the projections $\gamma$ and $\delta$ are 
always straight lines.  Any real curve is isometric with the graph of 
a real analytic function.

The three real curves which bound the triangular graph $P|Q$ generate 
a topological cell decomposition of $\Bar{P}|\Bar{Q}$ into vertices, 
edges, and the graph $P|Q$ itself.  The generators of the edge group, 
$\sigma_{ij}$, acting on $\Bar{P}|\Bar{Q}$, generate three closed graphs, each 
having one edge in common with $\Bar{P}|\Bar{Q}$.  By continuing this 
construction we obtain a cell decomposition of $[P|Q]$ into 
2-cells isometric with $P|Q$, 1-cells isometric with one of the edges 
of $\Bar{P}|\Bar{Q}$, and points.  The cell complex as a whole defines 
a tiling $\T$; the 2-cells by themselves form a set of tiles, 
$\T_{2}$, and every element of $\T_{2}$ can be expressed as $g(P|Q)$, 
where $g$ is an element of the edge group, $G$.

To show that $G$ is the maximal isometry group we first need to check 
that the tiling $\T$ is \emph{primitive}, 
that is, there is no refinement of the tiles $\T_{2}$ by 
additional real curves within $[P|Q]$ we may have missed.  For this it 
suffices to check that there are no real curves within $P|Q$.  Before 
we can prove this statement we need some basic properties of real 
curves.

\begin{defn}
A real curve is \emph{complete} if it is not a proper subset of any other 
real curve.
\end{defn}

\begin{lem}\label{singEndpoints}
The closure in $\Surf$ of a complete real curve 
$\gamma|\delta\subset\Check{\Surf}$, if bounded, has singular endpoints.
\end{lem}
\begin{proof}
Without loss of generality let
$\gamma$ and $\delta$ lie on the real axes of, respectively, 
$X$ and $Y$. The functions $f$ of the function elements $(U,f)$, 
which represent $\Surf$ 
locally, will then have power series on the real axis (of $X$) with 
real coefficients. Since a real power series when analytically 
continued along the real axis continues to be real, we can continue 
$\gamma|\delta$ until we encounter either a singularity of $f$ 
or a zero of $f'$ (\textit{i.e.} a singularity of $f^{-1}$ on the real axis of $Y$).
\end{proof}

The next Lemmas deal with the angles formed by intersecting real curves.
\begin{defn}
The \emph{angle between lines} $\gamma$ and $\gamma'$ (in $X$ 
or $Y$), 
denoted $\angle (\gamma, \gamma')$, is 
the smallest counterclockwise rotation required to make $\gamma$ 
parallel to $\gamma'$.
\end{defn}
\begin{lem}\label{sameAngle}
If real curves $\gamma|\delta\subset\Check{\Surf}$ and 
$\gamma'|\delta'\subset\Check{\Surf}$ intersect, then 
$\angle (\gamma, \gamma')=\angle (\delta, \delta')$.
\end{lem}
\begin{proof}
Near the point of intersection $\Check{\Surf}$ is represented by a 
function element $(U,f)$ where $f$ is conformal. The equality of 
angles, formed by a pair of lines in $X$ and their images by $f$ in 
$Y$, is simply the 
geometrical statement that $f$ is conformal.  
\end{proof}
\begin{lem}\label{intersectionAngle}
Only a finite number $n>1$ of real curves can intersect at any point 
of a nontrivial Riemann surface
and the angle formed by any pair must be a multiple of $\pi/n$.
\end{lem}
\begin{proof}
Suppose $\gamma|\delta$ and $\gamma'|\delta'$ intersect 
with angle $\angle (\gamma, \gamma')=\angle (\delta, 
\delta')=\alpha>0$ on a nontrivial Riemann surface $\Surf$;
for convenience, let $(0,0)$ be the point of intersection.
These curves are fixed by Schwarz
reflections $\sigma$ and $\sigma'$ respectively, and 
$\sigma'\sigma=r(2\alpha,2\alpha)$ is an isometry of $\Surf$. 
The neighborhood of the point of intersection is the graph
\begin{equation}
\Surf_{0}=\{(x,f(x))\colon x\in U\},
\end{equation}
where $U$ is a neighborhood of the origin in $X$, and $f$ is 
conformal at $x=0$. The Taylor series for $f$ at the origin has the 
form
\begin{equation}
f(x)=\sum_{k=1}^{\infty}a_{k}x^{k},
\end{equation}
where $a_{1}\neq 0$. A short calculation shows
\begin{equation}
r(2\alpha,2\alpha)\Surf_{0}=\{(x,f_{\alpha}(x))\colon x\in U_{\alpha}\},
\end{equation}
where $U_{\alpha}=e^{i 2\alpha}U$ is again a neighborhood of the 
origin, and
\begin{equation}
f_{\alpha}(x)=\sum_{k=1}^{\infty}a_{k}e^{i 2\alpha(1-k)}x^{k}.
\end{equation}
Since $r(2\alpha,2\alpha)$ is an isometry, the Taylor 
series for $f$ and $f_{\alpha}$ must agree, term by term. Now if 
$\alpha=\pi \omega$ and $\omega$ is irrational, then $\omega(1-k)$ 
can be an integer only for $k=1$ (so that $e^{i 2\alpha(1-k)}=1$). But this 
requires $a_{k}=0$ for $k>1$ which is impossible since $\Surf$ is 
nontrivial. Thus we may assume $\omega=p/q$ where $p$ and $q$ are 
relatively prime positive integers, $p<q$ (since $\alpha<\pi$), and 
$a_{1+m q}\neq 0$ for some integer $m>0$. 

Now let $\gamma''|\delta''$ 
be any real curve that intersects $\gamma|\delta$ at the origin; then
$\angle (\gamma, \gamma'')=\angle (\delta, 
\delta'')=\pi(p'/q')$ by the argument just given, where $p'$ and $q'$ are 
relatively prime positive integers, $p'<q'$. However, since $a_{1+m q}\neq 
0$, we must have $e^{-i 2\pi(p'/q')m q}=1$, or that $q'$ divides the 
product $m q=n$. This shows that $\angle (\gamma, \gamma'')$ is a 
multiple of $\pi/n$, for some $n>1$.
\end{proof}

Clearly any triangular graph with a nontrivial isometry must be 
``isoceles'' and fails to be primitive because it can be 
decomposed into two isometric tiles. This is made precise by the 
following Lemma.

\begin{lem}\label{noRealCurvesInside}
Let $P|Q$ be a nontrivial triangular graph with trivial isometry 
group, then $P|Q$ contains no real curves.
\end{lem}

\begin{figure}
\centerline{\epsfbox{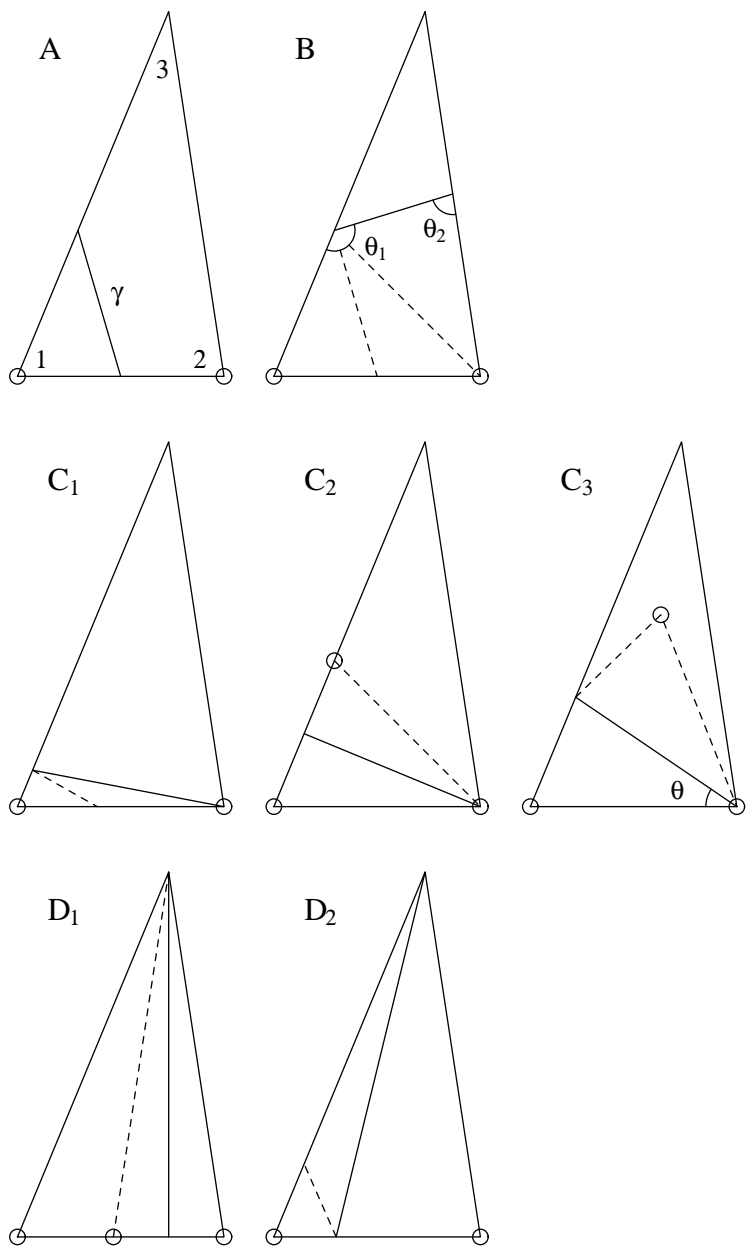}}
\caption{Diagrams used in the proof of Lemma \ref{noRealCurvesInside}.}
\end{figure}

\begin{proof}
Suppose $P|Q$ contains a real curve and call its completion 
$\gamma|\delta$. We recall that $\gamma$, and its closure in $X$, 
$\bar{\gamma}$, are straight lines and $\bar{\gamma}$
cannot have an endpoint within $P$ (Lemma \ref{singEndpoints}). 
The possible geometrical 
relationships between $\bar{\gamma}$ and $P$ are diagrammed in Figure 
2. 
Since $P|Q$ is nontrivial, at least two vertices are singular and are 
shown circled in each diagram. Either $\bar{\gamma}$ intersects two 
edges of $P$, as in cases $A$ and $B$, or, it intersects an edge and the  
opposite vertex which may be singular (case $C$) or possibly 
regular (case $D$). The vertex labels on the diagram refer to our 
notation for the vertex angles and edges. For example, $\alpha_{1}$ 
and $\beta_{1}$ are the angles in $P$ and $Q$, respectively, of 
vertex 1; $P_{12}|Q_{12}$ is the edge (real curve) bounded by vertices 
1 and 2, etc.

Case $A$ is easily disposed of using Lemma \ref{sameAngle}:
\begin{equation}
\begin{align}
\angle(\bar{\gamma},P_{12})=
\angle(\bar{\gamma},P_{13})+\alpha_{1}
=\angle(\bar{\delta},Q_{13})+\alpha_{1}
&=\angle(\bar{\delta},Q_{12})-\beta_{1}+\alpha_{1}\\
&=\angle(\bar{\gamma},P_{12})-\beta_{1}+\alpha_{1}.
\end{align}
\end{equation}
This is impossible because vertex 1 is singular 
($\alpha_{1}\neq\beta_{1}$).

By using Schwarz reflection to imply the existence of additional real 
curves, the remaining cases either reduce to case $A$ or imply the 
existence of a singularity within $P|Q$ or one of its edges --- 
neither of which is possible.

First consider case $B$. Let $\gamma$ intersect $P_{13}$ at 
$x_{1}$ and $P_{23}$ at $x_{2}$, forming angles $\theta_{1}$ and 
$\theta_{2}$ (see Fig. 2). Any other complete real curve with 
projection $\gamma'$ which intersects $x_{1}$ makes a finite angle 
with $\gamma$ by Lemma \ref{intersectionAngle}. Thus we may assume 
the angles $\theta_{1}$ and $\theta_{2}$ are the smallest possible 
(for a $\gamma$ that intersects 
both $P_{13}$ and $P_{23}$). Since one of $\theta_{1}$ and 
$\theta_{2}$ must be greater than $\pi/2$, we assume without loss of 
generality it is $\theta_{1}$. If we now reflect $\gamma$ in $P_{13}$ 
we obtain a real curve with projection $\gamma'$ such that $\gamma'$ 
intersects $P_{13}$ but not $P_{23}$. Thus case $B$ always reduces to 
cases $A$ or $C$.

In case $C$ we consider the sequence of real curves with 
projections $\gamma_{k}$, where $\gamma_{0}=P_{23}$, 
$\gamma_{1}=\gamma$, and $\gamma_{k+1}$ is the image of $\gamma_{k-1}$ 
under reflection in $\gamma_{k}$. Let $\theta_{k}$ be the angle formed 
at vertex 2 in $P$ by $\gamma_{k}$. Clearly for some $k$ 
we arrive at a $\gamma'=\gamma_{k}$ such that 
$\theta=\theta_{k}\leq\alpha_{2}/2$ (see Fig. 2). 
This leads to three subcases: $C_{1}$, where 
$\angle(\gamma',P_{13})\leq\pi/2$, $C_{2}$, where $\gamma'$ and 
$P_{13}$ are perpendicular, and $C_{3}$, where 
$\angle(\gamma',P_{13})\geq\pi/2$. In case $C_{2}$, $\theta<\alpha_{2}/2$ since 
otherwise $P|Q$ would have a nontrivial isometry (reflection in 
$\gamma'|\delta'$). All three subcases immediately lead to 
contradictions. In $C_{1}$, reflecting $\gamma'$ in $P_{13}$ presents us 
with a $\gamma''$ satisfying case $A$. In $C_{3}$, the image of vertex 1 
under reflection in $\gamma'$ implies a singularity within $P$; in 
$C_{2}$ 
the same reflection implies a singularity on $P_{13}$.

Case $D$: we either have $\angle(\gamma,P_{12})=\pi/2$, case $D_{1}$, or 
$\angle(\gamma,P_{12})\neq\pi/2$, case $D_{2}$. Since $P$ has no 
nontrivial isometry, a reflection in $\gamma$ in case $D_{1}$ would place the image 
of either vertex 1 or 2 (both singular) somewhere on $P_{12}$. In 
$D_{2}$, a 
Schwarz reflection 
of $\gamma$ leads to case $A$.
\end{proof}

\begin{thm}\label{maxIsometryGroup}
Let $P|Q$ be a nontrivial triangular graph with trivial isometry 
group, then the maximal isometry group of $[P|Q]$ is the edge group 
of $P|Q$. 
\end{thm}
\begin{proof}
We use the real curves to decompose $[P|Q]$ into 
a set of tiles (2-cells) $\T_{2}$. Lemma \ref{noRealCurvesInside} 
tells us that $P|Q\in\T_{2}$, so that $\T_{2}=G(P|Q)$, where $G$ is 
the edge group of $P|Q$. On the other hand,
if $h$ is an isometry of $[P|Q]$, then $h(P|Q)=P'|Q'\in\T_{2}$, 
where $P'|Q'=g(P|Q)$ for some $g\in G$. But since $P|Q$ has no 
nontrivial isometries, the map $h^{-1}g\colon P|Q\to P|Q$ must be the 
identity 
and $h=g$.
\end{proof}

We conclude this section with a formula for the topological genus of 
a Riemann surface $[P|Q]$ compactified by the translation subgroup of 
its isometry group, the \emph{lattice group} $\Lambda$ of $[P|Q]$.
\begin{defn}
The \emph{vertex groups} 
$G_{i}$, $(i=1,2,3)$, of a triangular graph $P|Q$, are the subgroups of the 
edge group of $P|Q$ generated by the adjacent edges of, respectively, 
the three vertices of 
$P|Q$.
\end{defn}
\begin{thm}
Let $P|Q$ be a nontrivial triangular graph with trivial isometry 
group. Let $G$ be the isometry group of $[P|Q]$, $\Lambda$ 
its lattice group, and $G_{i}$, $(i=1,2,3)$, the three vertex groups 
of $P|Q$. If $|G/\Lambda|$ is finite, the genus $g$ of the surface 
$[P|Q]/\Lambda$ satisfies
\begin{equation}\label{genusFormula}
2-2g=|G/\Lambda|\left(\sum_{i=1}^{3}\frac{1}{|G_{i}|}-\frac{1}{2}\right).
\end{equation}
\end{thm}
\begin{proof}
If $|G/\Lambda|$ is finite we can view $[P|Q]/\Lambda$ as a 
finite cell complex. We can relate the number of 0-cells, $N_{0}$, and 
the number of 1-cells, $N_{1}$, in this complex to the number of 
2-cells, $N_{2}$. Since every 2-cell is bounded by three 1-cells, each 
of which bounds exactly one other 2-cell, $N_{1}=(3/2)N_{2}$. 
Similarly, the boundary of each 2-cell contains three 
0-cells (the vertices $i=1,2,3$), each of which belongs to the 
boundary of a number of 2-cells equal to the order of the corresponding 
vertex group, $|G_{i}|$.  Thus 
$N_{0}=(\sum_{i=1}^{3}|G_{i}|^{-1})N_{2}$.  Finally, since $P|Q$ has 
no nontrivial isometry, and $G/\Lambda$ acts transitively on the 2-cells 
of
$[P|Q]/\Lambda$, $N_{2}=|G/\Lambda|$.  The result 
\eqref{genusFormula} follows from Euler's formula, 
$2-2g=N_{2}-N_{1}+N_{0}$.
\end{proof}

\subsection{Discreteness and uniformity}

The whole point of immersing a Riemann surface $\Surf$ in 
$X\times Y=\mathbb{C}^{2}$ is 
that by forming sections of $\Surf$, \textit{i.e.} intersections with $X=\mathbb{C}$, 
one obtains patterns of points. A very primitive 
property of a point set, normally taken for granted in 
crystallography, is discreteness.

\begin{defn}
The \emph{section} of the Riemann surface $\Surf$ at $x$ is the set
\begin{equation}
\A(x)=\pi_{Y}\circ {\pi_{X}}^{-1}(x).
\end{equation}
\end{defn}

The basic property of a holomorphic function, that its zeros form a 
discrete set, translates to the statement that the preimages 
${\pi_{X}}^{-1}(x)$ are discrete in a Riemann surface $\Surf$. We use the stronger 
property that ${\pi_{X}}^{-1}(x)$ is 
discrete in $X\times Y$ to define a \emph{discrete Riemann surface}.  
This is equivalent to the following statement about sections:

\begin{defn}
A Riemann surface $\Surf$ is \emph{discrete} if its sections $\A(x)$ 
are discrete in 
$Y$ for every $x\in X$.
\end{defn}

All the statements we can make about discreteness of a Riemann surface 
$\Surf$ hinge upon properties of the lattice group of $\Surf$, $\Lambda$. 
One property is the rank, $\Rk(\Lambda)$, given by 
the cardinality of the generators of $\Lambda$. 
The orbit of the origin of $X\times Y$, $\Lambda 
(0,0)$, is called a lattice and is also represented by the symbol 
$\Lambda$. When $\Rk(\Lambda)=4$, a second property is the
determinant of the lattice, $\det{\Lambda}$. If $\det{\Lambda}>0$, the 
four generators of $\Lambda$ are linearly independent (as vectors in 
$X\times Y$); if $\det{\Lambda}=0$ the generators are linearly 
dependent and $\Lambda$ (as a lattice) is not 
discrete in $X\times Y$. A lattice with $\Rk(\Lambda)>4$ is never 
discrete in $X\times Y$.

\begin{notation}
The standard measure for a set $A$ is written $|A|$. If $A$ is a 
region in $\mathbb{C}$ then $|A|$ is its area; if $A$ is a set of points, 
then $|A|$ is its cardinality.  Finally, if $\Lambda$ is a rank 4 
lattice, then $|\Lambda|=\sqrt{\det{\Lambda}}$ is the volume in 
$X\times Y$ of its fundamental region.
\end{notation}

The following Lemma provides a necessary condition for
discreteness:
\begin{lem}\label{notDiscrete}
The lattice group $\Lambda$, of a discrete, nontrivial Riemann 
surface $\Surf$, is discrete (as a lattice) in $X\times Y$ and in 
particular, $\Rk{\Lambda}\leq 4$.
\end{lem}
\begin{proof}
If $\Lambda$ is not discrete we can find a sequence of
$t(u,v)\in\Lambda$ such that both $\|u\|\to 0$ and $\|v\|\to 0$.
Let $(x_{0},y_{0})$ be a regular point of $\Surf$, then $y_{0}\in \A(x_{0})$.
Near $(x_{0},y_{0})$ we can represent 
$\Surf$ by the graph
\begin{equation}
U|V=\{(x,f(x))\colon x\in U\},
\end{equation}
where $U$ is a neighborhood of $x_{0}$ in $X$, $f$ is conformal in 
$U$, and 
$f(x_{0})=y_{0}$. Since $t(u,v)$ is an isometry, $t(u,v)(U|V)\subset 
\Surf$. From
\begin{equation}
t(u,v)(U|V)=\{(x+u,f(x)+v)\colon x\in U\},
\end{equation}
we see that $f(x_{0}-u)+v\in \A(x_{0})$, since $x_{0}-u\in U$ as 
$\|u\|$ can be arbitrarily small. 
Since $\Surf$ is discrete, $y_{0}$ is isolated in $Y$ and there must be a subsequence 
$t(u',v')$ such that $f(x_{0}-u')+v'=f(x_{0})$. If, within the 
sequence $t(u',v')$, there is a subsequence $t(u'',v'')$ with $u''=0$, 
then $v''=f(x_{0})-f(x_{0}-u'')=0$ and
we have a contradiction. Thus there must be a subsequence with 
$u''\neq 0$. Since $f$ is conformal at $x_{0}$,
\begin{equation}
\begin{align}
\lim_{u''\to 0}\frac{f(x_{0})-f(x_{0}-u'')}{u''}&=f'(x_{0})\\
&=\lim_{(u'',v'')\to (0,0)}\frac{v''}{u''}.
\end{align}
\end{equation}
But the second limit, above, is independent of $x_{0}$ so we are 
forced to conclude that $f'$ is constant. This is 
impossible because $\Surf$ is nontrivial.
\end{proof}

The point sets studied in crystallography normally are Delone sets 
and have the property of being \emph{uniformly discrete} \cite{Moody}.  
For a point set in $\mathbb{R}^{n}$ this means there exists a real 
number $r>0$ such that a spherical neighborhood of radius $r$ about 
any point of the set contains no other point of the set.  For the sets 
$\A(x)$ generated by Riemann surfaces this property is clearly too 
strong: it is violated whenever $x$ is near a branch point of 
$\pi_{X}$.  We therefore adopt a weaker form of this property which is 
nevertheless stronger than discreteness and useful in establishing the 
existence of the density.
\begin{defn}
Let $B_{r}(y)\subset Y$ be an open disk of radius $r$ centered at $y$. 
A Riemann surface $\Surf$, with sections $\A(x)$, is \emph{finitely discrete} 
if for some $r>0$, $|\A(x)\cap B_{r}(y)|$ is uniformly bounded above for 
all $x\in X$ and $y\in Y$.
\end{defn}
Since a disk of radius $r'$ can always be covered by finitely many 
disks of radius $r$, the finitely discrete property holds for any 
$r'$ once it has been established for a particular $r$. 

We now introduce the class of Riemann surfaces which is the focus of 
this study.
\begin{defn}
A Riemann surface
is \emph{crystallographic} if its lattice group $\Lambda$ has rank 4 and 
$|\Lambda|>0$.
\end{defn}

\begin{defn}
Any isometry $g$ of $X\times Y$ can uniquely be expressed in the form 
$g=g_{0}\lambda$, where $g_{0}$ fixes the origin and $\lambda$ is a 
translation.  The \emph{derived point group} of the isometry group $G$, 
$\psi(G)$, is the image of $G$ by the homomorphism $\psi\colon g\mapsto 
g_{0}$.
Since $\Ker{\psi}=\Lambda$, the lattice group of $G$, we have the 
isomorphism $\psi(G)\simeq G/\Lambda$.
\end{defn}
\noindent
It is important to remember that $\psi(G)$ need not be a subgroup 
of $G$; nevertheless, the lattice group of $G$ is always left 
invariant by $\psi(G)$, just as it is invariant within $G$.  
Furthermore, if a Riemann surface is crystallographic, then the action 
(by conjugation) of $\psi(G)$ on its lattice group $\Lambda$ is a faithful 
representation of $G/\Lambda$.  Since the isometry group of a (finite 
rank) lattice is finite, we have that $G/\Lambda$, for a 
crystallographic Riemann surface, is always finite.

\begin{lem}\label{uniformlyDiscrete}
A Riemann surface determined by a triangular graph $P|Q$, if crystallographic, 
is finitely discrete.
\end{lem}
\begin{proof}
Let $B_{r}(x,y)\subset X\times Y$ be an open ball of radius $r$ 
centered at an arbitrary point $(x,y)$. Consider the piece of the 
Riemann surface within this ball, 
$\Surf_{B}=[P|Q]\cap B_{r}(x,y)$, and the projection $\pi_{X}\colon 
\Surf_{B}\to X$. $[P|Q]$ is finitely discrete if there is a 
uniform upper bound on the number of preimages $\pi_{X}^{-1}(x)$.

$[P|Q]$ is covered by the orbit of closed graphs, 
$G(\Bar{P}|\Bar{Q})$, where $G$ is the edge group of $P|Q$.  Let 
$\Lambda$ be the lattice group of $G$, then $G(\Bar{P}|\Bar{Q})$ is the union 
of cosets, $H_{i}(\Bar{P}|\Bar{Q})$, $i=1,\ldots,N$, where 
$N=|G/\Lambda|$ is finite because $[P|Q]$ is crystallographic.  Again, 
because $[P|Q]$ is crystallographic, all but finitely many graphs in 
$H_{i}(\Bar{P}|\Bar{Q})$ have empty intersection with a ball of radius 
$r$, in particular, $B_{r}(x,y)$.  Thus we have a bound (independent of 
$x$ and $y$) on the number of graphs in $G(\Bar{P}|\Bar{Q})$ which 
intersect $B_{r}(x,y)$.  But a graph can have at most one preimage of 
$\pi_{X}$; hence $|\pi_{X}^{-1}(x)|$ is uniformly bounded above.
\end{proof}

The sections $\A(x)$ of a crystallographic Riemann surface also possess a 
uniformity with 
respect to the parameter $x$. Our handle on this property is provided, 
in part, by the smooth behavior of $\A(x)$ with $x$.  Before we can 
proceed, however, we need to be aware of two point sets in $X$ which 
create problems: branch points (of $\pi_{X}$) and crossing points.
\begin{defn}
A point $(x,y)\in\Surf$ is a \emph{self-intersection point} if in the 
description of $\Surf$ as a complete global analytic function there exist 
function elements $(U,f)$ and $(V,g)$, such that $x\in U\cap V$, 
$f\neq g$ in $U\cap V$, and $f(x)=g(x)=y$.  The point $x\in X$ is 
called a \emph{crossing point}.
\begin{lem}
A Riemann surface determined by a triangular graph $P|Q$ has countably many 
self-intersection points and $\pi_{X}$ has countably many branch 
points.
\end{lem}
\begin{proof}
$[P|Q]$ is the union of countably many closed graphs 
$\Bar{P}_{i}|\Bar{Q}_{i}$ given by the orbit of $\Bar{P}|\Bar{Q}$ 
under the action of the edge group. Since each graph has at most three 
singular points, $\pi_{X}$ has countably many branch points. If there 
were uncountably many self-intersection points then uncountably many 
must arise from one pair of distinct graphs, say $\Bar{P}_{i}|\Bar{Q}_{i}$ and 
$\Bar{P}_{j}|\Bar{Q}_{j}$. Let $f_{i}$ and $f_{j}$ be the 
corresponding conformal maps; then $f_{i}(x)=f_{j}(x)$ would have 
uncountably many solutions $x\in \Bar{P}_{i}\cap\Bar{P}_{j}$. 
Thus either $f_{i}=f_{j}$, a contradiction,
or the zeroes of $f_{i}-f_{j}$ would not be isolated, another 
impossibility.
\end{proof}
\end{defn}
\begin{lem}\label{uniformInX}
Let $[P|Q]$ be crystallographic, $X_{c}\subset X$ its crossing 
points, and $X_{b}\subset X$ the 
branch points of $\pi_{X}$; then for any pair $x,x'\in X\setminus 
(X_{b}\cup X_{c})$, there exists a bijection of sections of $[P|Q]$,
\begin{equation}
\Psi\colon \A(x)\to\A(x'),
\end{equation}
such that $\|y-\Psi(y)\|$ is uniformly bounded above for $y\in \A(x)$.
\end{lem}
\begin{proof}
We arrive at $\Psi$ by composing bijections
\begin{equation}
\begin{align}
\Psi_{1}&\colon\A(x)\to\A(x''),\\ 
\Psi_{2}&\colon\A(x'')\to\A(x'),
\end{align}
\end{equation}
such that 
$\|y-\Psi_{1}(y)\|$ and $\|y''-\Psi_{2}(y'')\|$ are (correspondingly) 
uniformly bounded. The Lemma then 
follows by application of the triangle inequality.

Since $[P|Q]$ is crystallographic, we can partition $X\times Y$ 
into translates of a bounded fundamental region, $V(0)$, of its 
lattice $\Lambda$. Thus for any pair $x,x'\in X$ we can write
\begin{equation}
\begin{align}
(x,0)&\in V(0)+\lambda,\\
(x',0)&\in V(0)+\lambda' ,
\end{align}
\end{equation}
where $\lambda,\lambda'\in\Lambda$. Consider the point
\begin{equation}\label{xDoublePrime}
(x'',y'')=(x,0)+\lambda'-\lambda\in V(0)+\lambda'.
\end{equation}
Since
\begin{equation}
(x''-x',y'')=(x'',y'')-(x',0)\in V(0)-V(0),
\end{equation}
both $\|x''-x'\|$ and $\|y''\|$ have upper bounds independent of $x$ and $x'$.
Now
\begin{equation}
\A(x)=\{y\in Y\colon (x,y)\in [P|Q]\},
\end{equation}
and
\begin{equation}
\A(x)+y''=\{y'\in Y\colon (x,y'-y'')\in [P|Q]\}.
\end{equation}
But $(x,y'-y'')\in [P|Q]$ iff $(x,y'-y'')+\lambda''\in [P|Q]$, where 
$\lambda''\in \Lambda$.  Choosing
\begin{equation}
\lambda''=\lambda'-\lambda=(x''-x,y''),
\end{equation}
we obtain
\begin{equation}
\A(x)+y''=\{y'\in Y\colon (x'',y')\in [P|Q]\}=\A(x'').
\end{equation}
As our first bijection we take the translation $\Psi_{1}(y)=y+y''$, 
where $\|y''\|$ is uniformly bounded from above. 
The point of this intermediate step is that for $\Psi_{2}$ we need 
consider only pairs of sections with bounded separation $\|x''-x'\|$.

In constructing $\Psi_{2}$ we avoid branch points and crossing points. 
Since $x\in 
X\setminus (X_{b}\cup X_{c})$, equation \eqref{xDoublePrime} implies $x''\in 
X\setminus (X_{b}\cup X_{c})$. Let $\gamma\colon [0,1]\to X\setminus (X_{b}\cup X_{c})$ 
be a smooth rectifiable curve with $\gamma(0)=x''$ 
and $\gamma(1)=x'$. To show that $\gamma$ exists we recall that 
$X_{b}\cup X_{c}$ is countable.
We can then find $\gamma$ in the uncountable family of 
circular arcs with endpoints $x''$ and $x'$, since each point of 
$X_{b}\cup X_{c}$ can eliminate at most one arc.

The curve $\gamma(t)$ generates a homotopy of the sections 
$\A(x'')$ and $\A(x')$. At each point $y''\in \A(x'')$, 
$\gamma(t)$ is lifted to a unique curve 
$\gamma(t)|\delta(t)\subset[P|Q]$ with endpoint 
$(\gamma(0),\delta(0))=(x'',y'')$ and we define our second 
bijection by $\Psi_{2}(y'')=\delta(1)$.  To finish the proof we need 
to show that $\|\delta(1)-\delta(0)\|$ is uniformly bounded.

Let $\Check{P}$ be the closed subset of $\Bar{P}$ that is a 
suitably small distance $r$ or greater from any of its vertices that are 
branch points of $\pi_{X}$. Let $\Check{P}|\Check{Q}\subset \Bar{P}|\Bar{Q}$
be the corresponding graph. The orbit under the edge group, 
$\Check{\Surf}=G(\Check{P}|\Check{Q})\subset 
[P|Q]$, is a Riemann surface from which all the points of 
ramification (of the map $\pi_{X}$) have been ``cut out''. The 
complement, $\Hat{\Surf}=[P|Q]\setminus\Check{\Surf}$, is the disjoint union of the 
branched neighborhoods of all the points of ramification.  It is 
possible to find curves $\gamma(t)$, such as the circular arcs 
considered above, where 
the branched neighborhoods $\Hat{\Surf}_{i}$ visited by 
$\gamma(t)|\delta(t)$ are visited only once, for 
$t\in T_{i}\subset [0,1]$.  Also, because we can bound the length $L$ 
of $\gamma$ and there is a minimum distance between branch points (on 
the branched covering of $X$), the number of such subintervals $T_{i}$ 
is bounded, \textit{i.e.} $i=1,\dots,N$.  The bounds on $N$ and $L$ 
are uniform bounds, independent of the points $x$ and $x'$.

Two additional bounds are needed before we can proceed to bound 
$\|\delta(1)-\delta(0)\|$. The first is an upper bound $D$ on 
the diameter of the projection of a 
branched neighborhood, $\pi_{Y}(\Hat{\Surf}_{i})$. This follows from 
the fact that $\Hat{\Surf}_{i}$ is isometric with the branched neighborhood 
$\Hat{\Surf}_{j}$ of a vertex of $\Bar{P}|\Bar{Q}$, and 
$\Hat{\Surf}_{j}\subset G_{j}(P|Q)$, where $G_{j}$ is the 
corresponding vertex group.  Clearly the maximum diameter of 
$\pi_{Y}(G_{j}(P|Q))$ is bounded because $Q$ is bounded.

The map $f\colon P\to Q$ (which defines 
$P|Q$), when restricted to $\Check{P}$ is conformal 
and $\|f'\|$ has a maximum value, $\mu$, since $\Check{P}$ is closed. 
This means that if $\gamma(t)|\delta(t)\in \Check{P}|\Check{Q}$, 
then
\begin{equation}
\left\|\frac{d\delta}{dt}\right\|=\left\|f'(\gamma)\frac{d\gamma}{dt}
\right\|\leq\mu\left\|\frac{d\gamma}{dt}\right\|.
\end{equation}
Because $\Check{\Surf}$ is generated from $\Check{P}|\Check{Q}$ by the 
action of $G$, this bounds applies globally, for $\gamma(t)|\delta(t)\in 
\Check{\Surf}$.

We are now ready to complete the proof:
\begin{equation}\label{mainInequality}
\|\delta(1)-\delta(0)\|=
\left\|\int_{[0,1]}\frac{d\delta}{dt}\right\|
\leq\left\|\int_{\cup_{i=1}^{N}T_{i}}\frac{d\delta}{dt}\right\|+
\left\|\int_{[0,1]\setminus\cup_{i=1}^{N}T_{i}}\frac{d\delta}{dt}\right\|.
\end{equation}
For each piece of the curve in a branched neighborhood we have
\begin{equation}
\left\|\int_{T_{i}}\frac{d\delta}{dt}\right\|\leq D,
\end{equation}
while in the complement ($\Check{\Surf}$),
\begin{equation}
\begin{align}
\left\|\int_{[0,1]\setminus\cup_{i=1}^{N}T_{i}}\frac{d\delta}{dt}\right\|&\leq
\int_{[0,1]\setminus\cup_{i=1}^{N}T_{i}}\left\|\frac{d\delta}{dt}\right\|\\
&\leq\mu\int_{[0,1]\setminus\cup_{i=1}^{N}T_{i}}\left\|\frac{d\gamma}{dt}\right\|\\
&\leq \mu L.
\end{align}
\end{equation}
Inequality \eqref{mainInequality} thus becomes
\begin{equation}
\|\delta(1)-\delta(0)\|\leq N D + \mu L.
\end{equation}
\end{proof}

For discrete Riemann surfaces with sufficiently uniform sections 
$\A(x)$, one can define their \emph{density}.
\begin{defn}
Let $B_{R}(0)\subset Y$ be a disk of radius $R$ centered at the origin. The limit
\begin{equation}
\rho(x)=\lim_{R\to \infty}\frac{|B_{R}(0)\cap \A(x)|}{|B_{R}(0)|},
\end{equation}
if it exists and is finite, is the density of $\A(x)$.
\end{defn}
With the aid of Lemmas \ref{uniformlyDiscrete} and \ref{uniformInX} 
we can show, that for a crystallographic Riemann surface generated by a triangular 
graph, 
$\rho(x)$ exists and is (essentially) independent of 
$x$.
\begin{notation}
The standard volume form in $X$ is $\omega_{X}=dx\wedge d\bar{x}$, its 
pullback on a Riemann surface $\Surf$ is written
${\pi_{X}}^{\ast}\omega_{X}$.
\end{notation}
\begin{thm}\label{density}
If $[P|Q]$ is crystallographic with lattice group $\Lambda$, 
its sections $\A(x)$ have density
\begin{equation}\label{KaluginFormula}
\rho=\frac{1}{|\Lambda|}\int_{[P|Q]/\Lambda}{\pi_{X}}^{\ast}\omega_{X},
\end{equation}
independent of $x$, provided $x$ is not a crossing point or a 
branch point of $\pi_{X}$. 
If $G$ is the edge group of $P|Q$, then
\begin{equation}\label{rhoFormula}
\rho=\frac{|G/\Lambda| |P|}{|\Lambda|}.
\end{equation}
\end{thm}
\begin{proof}
\begin{notation}
The expression $c=\order{(1/R)}$ indicates there exist constants $c_{1}$ 
and $c_{2}$ (independent of $R$) such that for sufficiently large $R$, 
$c_{1}/R<c<c_{2}/R$.
\end{notation}
Let $B_{R}(0)\subset Y$ be a disk of radius $R$ centered at the origin 
and let
\begin{equation}
N_{R}(x)=|B_{R}(0)\cap \A(x)|.
\end{equation}
We first obtain a bound on the difference,
$N_{R}(x)-N_{R}(x')$, when neither $x$ nor $x'$ 
is a crossing point or a branch point of $\pi_{X}$. 
By Lemma \ref{uniformInX} there exists 
a bijection $\Psi\colon \A(x)\to\A(x')$ such that if $y\in B_{R}(0)\cap 
\A(x)$, then $\Psi(y)\in B_{R+d}(0)\cap\A(x')$, where $d>0$ is a 
constant independent of $R$, $x$, and $x'$. This shows
\begin{equation}
\begin{align}
N_{R}(x)&\leq |B_{R+d}(0)\cap\A(x')|\\
&=N_{R}(x')+|(B_{R+d}(0)\setminus B_{R}(0))\cap\A(x')|.
\end{align}
\end{equation}
We can cover the annulus $B_{R+d}(0)\setminus B_{R}(0)$ by $M_{R}$ 
disks $B_{r}(y')$ of a fixed radius $r>0$, where, for sufficiently 
large $R$, $M_{R}<m R$ and $m$ is a constant independent of $R$. By 
Lemma \ref{uniformlyDiscrete}, $|B_{r}(y')\cap\A(x')|<n$, where $n$ is 
independent of $x'$ and $y'$. Thus $N_{R}(x)-N_{R}(x')<m n R$. 
Combining this bound with the bound obtained by 
interchanging $x$ and $x'$, we arrive at the statement
\begin{equation}\label{boundDeltaN}
N_{R}(x)-N_{R}(x')=|B_{R}(0)|\order{(1/R)}.
\end{equation}

We now introduce a disk $C_{R}(0)\subset X$ and consider the region 
$W_{R}=C_{R}(0)\times B_{R}(0)\subset X\times Y$. Since the set of 
branch points and crossing points is countable and has 
zero measure in $X$, and $\pi_{X}$ is otherwise smooth,
\begin{equation}\label{forms1}
\int_{\pi_{X}([P|Q]\cap W_{R})}\omega_{X}=
\int_{[P|Q]\cap W_{R}}{\pi_{X}}^{\ast}\omega_{X}.
\end{equation}
Because $[P|Q]$ is crystallographic, we can partition $X\times 
Y$ into translates of a bounded fundamental region of its lattice, 
$V(0)$. 
Let $\pi_{\Lambda}\colon [P|Q]\to [P|Q]/\Lambda$ be the standard 
projection on the quotient. On $[P|Q]\cap V(0)$ the map 
$\pi_{\Lambda}$ is 1-to-1 and
\begin{equation}\label{forms2}
\int_{[P|Q]\cap V(0)}{\pi_{X}}^{\ast}\omega_{X}=
\int_{[P|Q]/\Lambda}(\pi_{\Lambda}^{-1})^{\ast}{\pi_{X}}^{\ast}\omega_{X}=\rho |\Lambda|.
\end{equation}
This defines $\rho$, which we can make positive by appropriate choice 
of orientation on $[P|Q]$.

Turning now to the region $W_{R}$, there is a maximal subset 
$\Lambda_{-}\subset \Lambda$ such that $\Lambda_{-}+V(0)\subset W_{R}$ and a 
smallest subset $\Lambda_{+}\subset \Lambda$ such that 
$W_{R}\subset\Lambda_{+}+V(0)$. If $\lambda\in\Lambda_{+}\setminus\Lambda_{-}$, 
then $V_{\lambda}=(\lambda+V(0))\cap W_{R}$ is a proper subset of a 
fundamental region and
\begin{equation}
0<\int_{[P|Q]\cap 
V_{\lambda}}{\pi_{X}}^{\ast}\omega_{X}<\rho|\Lambda|.
\end{equation}
From this it follows that
\begin{equation}
\rho|\Lambda_{-}||\Lambda|<
\int_{[P|Q]\cap W_{R}}{\pi_{X}}^{\ast}\omega_{X}<
\rho|\Lambda_{+}||\Lambda|,
\end{equation}
and, from straightforward estimates of $\Lambda_{+}$ and $\Lambda_{-}$, we 
conclude
\begin{equation}\label{forms3}
\int_{[P|Q]\cap 
W_{R}}{\pi_{X}}^{\ast}\omega_{X}=\rho|W_{R}|(1+\order{(1/R)}).
\end{equation}

The projection $\pi_{X}([P|Q]\cap W_{R})$ covers the disk 
$C_{R}(0)$ multiple times, the multiplicity at the point $x\in 
C_{R}(0)$ being the number $N_{R}(x)$ defined above. Thus
\begin{equation}\label{forms4}
\int_{\pi_{X}([P|Q]\cap W_{R})}\omega_{X}=
\int_{C_{R}(0)}N_{R}(x)\omega_{X}.
\end{equation}
We can again neglect the countable set of branch points $X_{b}$ and 
crossing points $X_{c}$ to argue, 
for $x_{0}\in X\setminus (X_{b}\cup X_{c})$ fixed,
\begin{equation}\label{forms5}
\begin{align}
\int_{C_{R}(0)}N_{R}(x)\omega_{X}&=
\int_{C_{R}(0)}N_{R}(x_{0})\omega_{X}+
\int_{C_{R}(0)}(N_{R}(x)-N_{R}(x_{0}))\omega_{X}\nonumber\\
&=N_{R}(x_{0})|C_{R}(0)|+|C_{R}(0)||B_{R}(0)|\order{(1/R)},
\end{align}
\end{equation}
where in the last step we used \eqref{boundDeltaN}. Combining 
\eqref{forms1}, \eqref{forms3}, \eqref{forms4}, \eqref{forms5}, and 
using $|W_{R}|=|B_{R}(0)||C_{R}(0)|$, we obtain
\begin{equation}
\frac{N_{R}(x_{0})}{|B_{R}(0)|}=\rho(1+\order{(1/R)}),
\end{equation}
and thus
\begin{equation}
\lim_{R\to\infty}\frac{N_{R}(x_{0})}{|B_{R}(0)|}=\rho.
\end{equation}
To evaluate $\rho$ from \eqref{forms2}, we regard $[P|Q]/\Lambda$ 
as $|G/\Lambda|$ equivalence classes of tiles, all isometric to $P|Q$. 
The result \eqref{rhoFormula} follows because the integral of the 
form ${\pi_{X}}^{\ast}\omega_{X}$ over $P|Q$ is just the volume of 
$\pi_{X}(P|Q)=P$ in $X$.
\end{proof}

Formula \eqref{KaluginFormula} for the density was introduced by 
Kalugin \cite{Kalugin} to extend the notion of stoichiometry to quasicrystals.  
Because the 2-form ${\pi_{X}}^{\ast}\omega_{X}$ is closed, this 
formula gives the same density for 2-manifolds homologous in the torus 
$(X\times Y)/\Lambda$.  One must remember, however, that this homology 
invariant only corresponds to the true density when the map $\pi_{X}$ 
is orientation preserving (see equation \eqref{forms4}).

\newpage

\section{Classification of discrete Riemann surfaces generated by 
conformal maps of right triangles}

\subsection{Conformal maps of right triangles}

The simplest nontrivial conformal maps of triangles, $f\colon P\to Q$,  
are those where one of 
the vertices of the corresponding graph $P|Q$ is regular. The
edges of $P|Q$ adjacent to this vertex are real 
curves, and by Lemma \ref{intersectionAngle}, belong to a set of $n>1$ real 
curves intersecting at the same vertex with minimum angle $\pi/n$. 
This suggests that among those maps with one regular vertex, 
the simplest case is $n=2$, \textit{i.e.} the conformal map of 
right triangles with the right angle being the regular vertex.

Without loss of 
generality, we give $P$ and $Q$ a standard scale, position, and 
angular orientation as specified by the vertices $x_{i}|y_{i}$ of the 
corresponding graph $P|Q$:
\begin{equation}\label{triangleVertices}
\begin{align}
x_{1}|y_{1}&=0|0\nonumber\\
x_{2}|y_{2}&=\cos{\alpha}|\cos{\beta}\\
x_{3}|y_{3}&=e^{i\alpha}|e^{i\beta},\nonumber
\end{align}
\end{equation}
where $0<\alpha<\pi/2$ and 
$0<\beta<\pi/2$ are two free real parameters. Below we frequently use 
the abbreviations $a=\cos{\alpha}$, $b=\cos{\beta}$.
$P$ and $Q$ have angles $\alpha$ and $\beta$, respectively, at vertex 
1, and the corresponding complementary angles at vertex 3.
A nontrivial graph $P|Q$ has $\alpha\neq\beta$ with vertices 1 and 3 
singular; vertex 2 is regular.

The Schwarz-Christoffel formula \cite{Ahlfors} gives the conformal map 
$f$ as the composition $f=h\circ g^{-1}$ where $g$ and $h$ map the 
upper half plane of $Z=\mathbb{C}$ conformally onto, respectively, $P$ 
and $Q$.  Explicitly:
\begin{align}
g(z)&=A\int_{0}^{z}z^{\frac{\alpha}{\pi}-1}
(1-z)^{-\frac{1}{2}}dz,\label{map_{g}}\\
h(z)&=B\int_{0}^{z}z^{\frac{\beta}{\pi}-1}
(1-z)^{-\frac{1}{2}}dz.\label{map_{h}}
\end{align}
In both \eqref{map_{g}} and \eqref{map_{h}} the branches of the 
fractional powers are chosen so that the integrands are real and 
positive for 
$z\in (0,1)$. The normalization factors $A$ and $B$ are positive real 
numbers determined by the conditions $g(1)=a$ and $h(1)=b$. Further properties 
of $g$ and $h$ are easily checked, in 
particular, $g(\infty)=e^{i\alpha}$ and $h(\infty)=e^{i\beta}$. 

\subsection{Isometry groups}

\begin{notation}
Let $K$ be a set of elements of a group $G$, and let $k\in G$ be some 
element. We denote by $\langle K\rangle$ 
the subgroup generated by the elements of $K$, and by $\{k\}_{G}$ the 
conjugacy class of $k$ in $G$.
\end{notation}
In the case of right triangles, $P|Q$ can have a nontrivial isometry only if 
$\alpha=\pi/4$ and $\beta=\pi/4$.  But this makes $P|Q$ trivial.  From 
Theorem \ref{maxIsometryGroup} we know that for nontrivial $P|Q$ the 
isometry group $G$ is just the edge group generated by the three 
Schwarz reflections:
\begin{equation}
\begin{align}
\sigma_{12}&=\sigma,\\
\sigma_{13}&=r(2\alpha,2\beta) \sigma,\\
\sigma_{23}&=t(2 a,2 b) r(\pi,\pi) \sigma.
\end{align}
\end{equation}
From the isometry group
\begin{equation}\label{generators}
G=\langle\sigma_{12},\sigma_{13},\sigma_{23}\rangle
=\langle\sigma,r(2\alpha,2\beta),t(2 a,2 b)r(\pi,\pi)\rangle,
\end{equation}
we wish to extract the lattice group $\Lambda$. Helpful in this 
enterprise are the vertex group
\begin{equation}\label{G1}
G_{1}=\langle\sigma,r(2\alpha,2\beta)\rangle,
\end{equation}
and its cyclic subgroup,
\begin{equation}\label{R}
R=\langle r(2\alpha,2\beta)\rangle.
\end{equation}
Sets of translations invariant with respect to $G_{1}$, or \emph{stars}, play 
a central role in the construction of $\Lambda$. In what follows we 
will need two stars: 
\begin{equation}
\begin{align}
\Sigma&=\{t(2 a,2 b)\}_{G_{1}},\label{Sigma}\\
\Sigma^{-1}&=r(\pi,\pi)\,\Sigma\,r(\pi,\pi)=\{t(-2 a,-2 
b)\}_{G_{1}}.\label{Sigma-1}
\end{align}
\end{equation}
The main result is contained in the following Lemma:
\begin{lem}\label{groups}
Let $G$ be the isometry group of the Riemann surface $[P|Q]$ 
generated by the conformal map of the right triangles specified in 
\eqref{triangleVertices} and let $R$ be defined by \eqref{R}, $G_{1}$ 
by \eqref{G1}, $\Sigma$ and $\Sigma^{-1}$ by \eqref{Sigma} and 
\eqref{Sigma-1}.  If $r(\pi,\pi)\in 
R$, then $G$ has lattice group $\Lambda=\langle\Sigma\rangle$ and 
$G=\Lambda G_{1}$; otherwise, $G$ has lattice group 
$\Lambda=\langle\Sigma\Sigma^{-1}\rangle$ and $G=\Lambda G_{1}\, \cup\, t(2 
a,2 b) r(\pi,\pi)\Lambda G_{1}$.
\end{lem}
\begin{proof}
First consider the case $r(\pi,\pi)\in R$; then
\begin{equation}
G=\langle\sigma,r(2\alpha,2\beta),t(2 a,2 b)\rangle.
\end{equation} 
Now consider the group $H=\Lambda G_{1}$, where $\Lambda=\langle \Sigma\rangle$ 
is normal in $H$. Clearly $H\subset G$. Moreover, one easily verifies $gH=Hg=H$, 
where $g$ is any 
of the
three generators of $G$.  These two facts together show $G=H$; 
$\Lambda$ is clearly the lattice group of $G$.

Next consider the case $r(\pi,\pi)\notin R$. For the generators of 
$G$ we must now use \eqref{generators}. Consider the group
$\Tilde{G}=\Tilde{\Lambda}G_{1}$, where $\Tilde{\Lambda}=\langle \Sigma 
\Sigma^{-1}\rangle$ is normal in $\Tilde{G}$.  Clearly 
$\Tilde{G}\subset G$.  In contrast to the previous case, we can now 
only verify that $g\Tilde{G}g^{-1}=\Tilde{G}$, where $g$ is any of the 
three generators in \eqref{generators}.  Thus $\Tilde{G}$ is normal in 
$G$.  $\Tilde{G}$ has index at most two, since multiplication of 
$\Tilde{G}$ by the generators of $G$ produces at most two, possibly 
distinct, cosets: $\Tilde{G}$ and $\Tilde{G}'=t(2 a,2 b) 
r(\pi,\pi)\Tilde{G}$.  But if $\Tilde{G}'=\Tilde{G}$, then we would 
have some $\tilde{\lambda}\in\Tilde{\Lambda}$ and some $g_{1}\in G_{1}$ such that 
$\tilde{\lambda} g_{1}=t(2 a,2 b) r(\pi,\pi)$, or $g_{1}=\tilde{\lambda}^{-1} t(2 a,2 
b) r(\pi,\pi)$.  Since $g_{1}$ fixes the origin, $\tilde{\lambda}^{-1} t(2 a,2 
b)$ must be the trivial translation and $g_{1}=r(\pi,\pi)$.  This 
contradicts our assumption $r(\pi,\pi)\notin R$ and we conclude that 
$\Tilde{G}$ and $\Tilde{G}'$ are distinct.  Let $\Lambda$ be the 
lattice group of $G$.  Clearly $\Tilde{\Lambda}\subset\Lambda$. Now suppose 
$\lambda\in\Lambda$ but $\lambda\notin\Tilde{\Lambda}$. Since then 
$\lambda\notin\Tilde{G}$, we must have $\lambda\in\Tilde{G}'$, that 
is, $\lambda=t(2 a,2 b)r(\pi,\pi)g_{1}\tilde{\lambda}$ for some 
$g_{1}\in G_{1}$ and $\tilde{\lambda}\in\Tilde{\Lambda}$.  But this implies
$r(\pi,\pi)g_{1}=t(-2 a,-2 b)\lambda{\tilde{\lambda}}^{-1}$, a translation, 
and we arrive at the contradiction $r(\pi,\pi)g_{1}=1$.  Thus $\Lambda=\Tilde{\Lambda}$.
\end{proof}

\subsection{The discreteness restriction}

The requirement that $[P|Q]$ is discrete places strong constraints on 
the angles $\alpha$ and $\beta$ of the right triangles $P$ and $Q$.

\begin{lem}\label{rationalAngle}
If $[P|Q]$ is discrete, then $\alpha$ and $\beta$ are 
rational multiples of $\pi$.
\end{lem}
\begin{proof}
By Lemma \ref{notDiscrete}, $[P|Q]$ has a discrete lattice, \textit{i.e.} 
the orbit $\Lambda\cdot (0,0)$ is discrete in $X\times Y$. For right 
triangles, Lemma \ref{groups} gives us $\Lambda=\langle \Sigma\rangle$ 
if $r(\pi,\pi)\in R$, 
$\Lambda=\langle\Sigma\Sigma^{-1}\rangle$ otherwise. Thus 
discreteness of $[P|Q]$ implies discreteness of the star $\Sigma\cdot (0,0)$ 
in $X\times Y$. Since $\Sigma=\{t(2 a, 2 b)\}_{G_{1}}=\{t(2 a, 2 
b)\}_{R}$, $\Sigma\cdot (0,0)$ is just the orbit of $t(2 a,2 b)\cdot 
(0,0)$ under action of the group $R$ generated by $r(2\alpha,2\beta)$. 
Clearly $\Sigma\cdot (0,0)$ lies in a 2-torus $S_{1}\times S_{1}$ 
embedded in $X\times Y$. Since $\Sigma\cdot (0,0)$ is discrete, there 
is a disjoint union of neighborhoods, each containing just one element 
of $\Sigma\cdot (0,0)$. Moreover, since $R$ is an isometry of $X\times 
Y$ and acts transitively on $\Sigma\cdot (0,0)$, there is a uniform 
lower bound on the volumes of these neighborhoods. This implies the 
existence of disjoint neighborhoods in $S_{1}\times S_{1}$, again with 
a uniform lower bound on their measure. Since $S_{1}\times S_{1}$ has 
finite measure, this is only possible if $R$ has finite order.
\end{proof}

\begin{notation}
Given positive integers $m$ and $n$, $\GCD{(m,n)}$ is their greatest 
common divisor, $\LCM{(m,n)}$ their least common multiple.
\end{notation}
Lemma \ref{rationalAngle} allows us to write $\alpha=\pi (i/k)$, 
$\beta=\pi (j/l)$, where $i$, $j$, $k$ and $l$ are positive integers 
and $\GCD{(i,k)}=\GCD{(j,l)}=1$. Since $R=\langle r(2\pi(i/k),2\pi(j/l))\rangle$, 
we identify $n=\LCM{(k,l)}$ as the order of $R$. A more convenient 
parameterization is given by
\begin{equation}
\alpha=\pi\frac{p}{n},\quad\quad\beta=\pi\frac{q}{n},\label{angles}
\end{equation}
where $p=i(n/k)$, $q=j(n/l)$ are positive integers with no common factors that 
are also factors of $n$.

The transformation
\begin{equation}\label{complementaryAngles}
t(\sin{\alpha},\sin{\beta})\,r(-\pi/2,-\pi/2)\,t(-\cos{\alpha},-\cos{\beta})\,\sigma,
\end{equation} 
has the effect of replacing the 
angles $\alpha$ and $\beta$
by $\alpha'=\pi/2-\alpha$ and $\beta'=\pi/2-\beta$ in the definition 
of the triangles $P$ and $Q$.
Without loss of 
generality we may therefore take the smallest of the four 
angles $\alpha$, $\alpha'$, $\beta$ and $\beta'$, and rename this 
angle $\alpha$.  Transformation \eqref{complementaryAngles}, 
as well as an interchange of the spaces 
$X$ and $Y$, will then give all other cases of triangles $P$ and $Q$ 
(in their standard position, orientation and scale).
From $0<\alpha\leq \alpha'$ we obtain
\begin{equation}
0<p\leq\frac{n}{4}.\label{pInequality}
\end{equation}
The inequalities $\alpha\leq\beta$ and $\alpha\leq\beta'$, plus 
the condition to avoid triviality, $\alpha\neq\beta$, then give
\begin{equation}
p<q\leq\frac{n}{2}-p.\label{qInequality}
\end{equation}
Together, inequalities \eqref{pInequality} and \eqref{qInequality} 
have no solution unless
\begin{equation}
n\geq 6.\label{lowerBound}
\end{equation}

Since the vertex group $G_{1}$ is generated by $R$ and $\sigma$,
\begin{equation}\label{orderG1}
|G_{1}|=2 n.
\end{equation}
We can use transformation \eqref{complementaryAngles}, which has the 
effect of interchanging the angles at vertices 1 and 3, to compute the 
order of $G_{3}$.  Let $n'$, $p'$ and $q'$ be the integers 
parameterizing the angles at vertex 3 (as in \eqref{angles}); then
\begin{equation}\label{switch1and3}
n'=\frac{2 n}{\GCD{\left(2 n,n-2 p,n-2 q\right)}}
\end{equation}
\begin{equation}
\begin{align}
\frac{p}{n}+\frac{p'}{n'}&=\frac{1}{2}\nonumber\\
\frac{q}{n}+\frac{q'}{n'}&=\frac{1}{2},\nonumber
\end{align}
\end{equation}
and,
\begin{equation}\label{orderG3}
|G_{3}|=2 n'.
\end{equation}
Finally, since the group of the regular vertex
is generated by two real curves which intersect at right angles, 
\begin{equation}\label{orderG2}
|G_{2}|=4.
\end{equation}

From 
\eqref {generators} we find
\begin{equation}
\psi(G)=\langle G_{1}, r(\pi,\pi) \rangle=\langle \sigma, R, 
r(\pi,\pi) \rangle
\end{equation}
for the derived point group of $G$.
Recognizing $n=|R|$ as the order of an element of the 
isometry group of the \emph{lattice} $\Lambda$, we can use the 
following theorem of Senechal \cite{Senechal} and Hiller \cite{Hiller}, and the 
requirement of discreteness, to bound $n$ from above.
\begin{thm}[Senechal  \cite{Senechal}, Hiller \cite{Hiller}]\label{smallestRank}
Let $N(n)$ be the smallest integer such that the group 
$\GL{(N(n),\mathbb{Z})}$ has an element of order $n$, then
\begin{equation}
N(n)=\sum_{{p_{i}}^{m_{i}}\neq 2}\phi({p_{i}}^{m_{i}}),
\end{equation}
where ${p_{1}}^{m_{1}}\dotsm$ is the prime factorization of $n$ and
$\phi(k)$ is Euler's totient function: the number of positive 
integers less than and relatively prime to $k$ (note: the prime $2$ is 
included in the sum only if its exponent is greater than $1$).
\end{thm}
\begin{lem}\label{nValues}
If a Riemann surface $[P|Q]$ generated by the conformal map of 
right triangles $P$ and $Q$ is discrete, then the angles of $P$ and $Q$ are given 
by \eqref{angles} and either $n\leq 6$, or $n=8,10$ or $12$.
\end{lem}
\begin{proof}
Let $G$ be the isometry group of $[P|Q]$, $\psi(G)$ its derived 
point group, and $\Lambda$ its lattice group.
Since conjugation by $g\in\psi(G)$ leaves $\Lambda$ invariant, 
consider the automorphisms $\Phi_{g}\colon 
\Lambda\to\Lambda$ given by $\Phi_{g}(\lambda)=g\lambda g^{-1}$.  If 
$\Lambda$ has rank $N$, then the homomorphism $\Psi\colon\psi(G)\to 
\Aut(\Lambda)$, where $\Psi(g)=\Phi_{g}$, induces a representation of 
$\psi(G)$ by integral $N\times N$ matrices of determinant $\pm 1$.  In 
fact, $\Psi$ is an isomorphism since for any of the generators $g$ of 
$\psi(G)$ we can easily find a $\lambda\in \Lambda$ which is not fixed 
by $\Phi_{g}$ (it suffices to look within the star $\Sigma$ or, if 
$r(\pi,\pi)\notin R$, $\Sigma \Sigma^{-1}$).  Since 
$R\subset\psi(G)$ has order $n$, there must be an element of order $n$ 
in $\GL{(N,\mathbb{Z})}$.  By Theorem \ref{smallestRank} we must have 
$N\geq N(n)$.  On the other hand, if $[P|Q]$ is discrete, then 
$\Lambda$ must be discrete (as a lattice in $X\times Y$) which is 
possible only if $N\leq 4$.  Since $\phi(p^{m})=p^{m-1}(p-1)$, we need only 
consider the values of $\phi$ for powers of small primes: $\phi(4)=2$, 
$\phi(8)=4$, $\phi(3)=2$, $\phi(5)=4$ (all other primes and higher 
powers yield values greater than $4$).  From these facts we obtain 
just the values of $n$ given in the statement of the Lemma.
\end{proof}

With Lemma \ref{nValues} and inequalities \eqref{pInequality} and 
\eqref{qInequality}, the set of possible combinations of $n$, $p$ and 
$q$ is already finite.  Several of these combinations can be 
eliminated by the following Lemma which provides a lower bound on 
$\Rk(\Lambda)$ when either $p$ or $q$ is a nontrivial divisor of $n$.

\begin{lem}\label{rankBound}
Let $\Lambda$ be the lattice group of a discrete Riemann surface 
$[P|Q]$ with $P$ and $Q$ defined by parameters $n$, $p$ and $q$ 
\eqref{angles}; then if $p>1$ and $p$ divides $n$, 
$\Rk(\Lambda)\geq\phi(n/p)+\phi(n/\GCD{(n,q)})$ (and the 
same statement with $p$ and $q$ interchanged).
\end{lem}
\begin{proof}
Suppose $p$ divides $n$ and $n/p=d>1$. We may assume that $p$ does not divide 
$q$ since otherwise $n$, $p$ and $q$ would have $p$ as a common 
divisor. Let 
$r=r(2\alpha,2\beta)=r(2\pi/d,2\pi(q/n))$; then $R=\langle r\rangle$ 
and $R_{p}=\langle r^{d}\rangle$ is a subgroup of $R$ of order $p$. By 
looking in $\Sigma$ (or $\Sigma\Sigma^{-1}$), we can find a translation 
$t_{0}=t(u,v)\in\Lambda$ such that $u\neq 0$ and $v\neq 0$.  Now 
consider the two products of translations
\begin{equation}
\begin{align}
t_{1}&=t(u_{1},v_{1})=\prod_{s\in R_{p}}(s\:t_{0}\:s^{-1}),\\
t_{2}&=t(u_{2},v_{2})=\prod_{k=1}^{d}(s_{k}\:t_{0}\:{s_{k}}^{-1}),
\end{align}
\end{equation}
where $s_{k}$ is any element of the coset $r^{k}R_{p}$, and
$t_{2}$ depends on the particular choice of coset elements. 
Evaluating the products we find
\begin{equation}
\begin{align}
u_{1}&=p\,u,\\
v_{1}&=\left(\sum_{m=1}^{p}e^{2\pi i(q m/p)}\right)v.\label{zeroSum}
\end{align}
\end{equation}
Equation \eqref{zeroSum} implies $e^{-2\pi 
i(q/p)}v_{1}=v_{1}$, and, since $p$ does not divide $q$, we conclude 
$v_{1}=0$. Thus $t_{1}=t(p u,0)$. Similarly, we find
\begin{equation}
u_{2}=\left(\sum_{k=1}^{d}e^{2\pi i(k/d)}\right)u=0.
\end{equation}
On the other hand, $v_{2}$ is changed just by making a 
different choice for one coset element $s_{k}$ (again because $p$ 
does not divide $q$). Thus we can always make a choice such that 
$t_{2}=t(0,v_{2})$, where $v_{2}\neq 0$.

Now $\langle\{t_{1}\}_{R}\rangle=\Lambda_{X}$ is a lattice in $X$ isomorphic 
to the cyclotomic lattice $\mathbb{Z}[e^{2\pi i/d}]$, while 
$\langle\{t_{2}\}_{R}\rangle=\Lambda_{Y}$ is a lattice in $Y$ isomorphic to
$\mathbb{Z}[e^{2\pi i(q/n)}]$. Since 
$\Lambda\supset\Lambda_{X}\Lambda_{Y}$, 
$\Rk(\Lambda)\geq\Rk(\Lambda_{X})+\Rk(\Lambda_{Y})$. The statement of 
the Lemma follows from the well known formula for the rank of a 
cyclotomic lattice.
\end{proof}

As an example of the application of Lemma \ref{rankBound}, consider 
the case $n=8$, $p=1$ and $q=2$. 
For these numbers Lemma \ref{rankBound} gives 
$\Rk(\Lambda)\geq\phi(4)+\phi(8)=6$, and the corresponding Riemann 
surface would not be discrete.
Together, Lemma \ref{nValues}, Lemma \ref{rankBound} and the inequalities 
\eqref{pInequality} and \eqref{qInequality} limit the set of 
possibilities for $n$, $p$, and $q$ to the combinations 
$(n,p,q)=(6,1,2)$, $(8,1,3)$, $(10,1,3)$, $(10,1,4)$, $(12,1,5)$, and 
$(12,2,3)$.  All other combinations (consistent with the inequalities) 
correspond to Riemann surfaces whose lattice groups have ranks 
exceeding 4 and therefore cannot be discrete.  We have no general 
method to settle the discreteness of these remaining candidates and 
therefore considered them case by case.  The procedure was to obtain 
generators $e_{i}$ of the finite sets $\Sigma$, or $\Sigma\Sigma^{-1}$ 
if $r(\pi,\pi)\notin R$, since these generate the lattice group 
$\Lambda$.  The computations were performed in \textit{Mathematica} 
using the function \texttt{LatticeReduce}, which implements the 
Lenstra-Lenstra-Lovasz (LLL) lattice reduction algorithm.  With the 
exception of the case $(n,p,q)=(10,1,4)$, which was found to have rank 
8, all others proved to have rank 4 and nonvanishing determinant.  To 
these 5 cases of crystallographic Riemann surfaces, two more should be 
added which simply correspond to an interchange of the spaces $X$ and 
$Y$ (obtained by interchanging $p$ and $q$).  Although the interchangeability of these spaces was assumed in 
the derivation of inequalities \eqref{pInequality} and 
\eqref{qInequality}, clearly the sections $\A(x)$ will detect a 
difference. For example, there is no isometry which relates the 
surfaces $(10,1,3)$ and $(10,3,1)$ (or their sections). On the other 
hand, the use of transformation \eqref{complementaryAngles} shows that the 
surface $(10,3,1)$ is isometric with the surface described by 
$(n',p',q')=(5,1,2)$ (where now $p'<q'$). However, not all 
interchanges of $p$ and $q$ produce a 
new Riemann surface.  For example, $(n,p,q)=(8,3,1)$ 
corresponds, by \eqref{switch1and3}, to $(n',p',q')=(8,1,3)$ (the 
values before the interchange).

We summarize this discussion by our main Theorem:
\begin{thm}
Up to linear transformations, there are seven discrete Riemann surfaces
$[P|Q]$ generated by 
conformal maps of right triangles $P$ and $Q$. Each of these surfaces 
is crystallographic; their properties, in particular the integers 
$n$, $p$ and $q$ which specify the angles of $P$ and $Q$, are given
in Table 1.
\end{thm}
\begin{proof}
The properties listed in Table 1 are simple consequences  
of general results.
Since 
\begin{equation}
G/\Lambda\simeq\psi(G)=\langle\sigma, R, 
r(\pi,\pi)\rangle,
\end{equation}
we see that $G/\Lambda$ is isomorphic to an abstract group generated by 
an element of order 2, $\tilde{\sigma}$, and an element $\tilde{r}$, 
which has order $2|R|=2 n$, if $n=|R|$ is odd (and therefore 
$r(\pi,\pi)\notin R$) or $n$ is even and just one of $p$ and $q$ is 
odd (since then the element of order two in $R$ is either $r(\pi,0)$ 
or $r(0,\pi)$).  Otherwise ($n$ even, $p$ and $q$ both odd), 
$r(\pi,\pi)\in R$ and $\tilde{r}$ has order $|R|=n$. These generators have the 
relation $\tilde{\sigma}\tilde{r}\tilde{\sigma}={\tilde{r}}^{-1}$ and 
imply $G/\Lambda\simeq d_{n}$, the dihedral group of order $2 n$, 
or $G/\Lambda\simeq d_{2 n}$.

To compute the genus we use the orders of the vertex groups, 
\eqref{orderG1}, \eqref{orderG3} and \eqref{orderG2}, in formula \eqref{genusFormula}:
\begin{equation}
2-2g=|G/\Lambda|\left(\frac{1}{2n}+\frac{1}{2n'}-\frac{1}{4}\right).
\end{equation}

The geometry of $\Lambda$ 
is completely specified by the Gram matrix formed from its generators 
$e_{k}$, $k=1,\ldots,4$:
\begin{equation}
(M_{kl})=e_{k}\cdot e_{l},
\end{equation}
where $\cdot$ is the standard inner product.
If we consider each $e_{k}$ as a vector in $\mathbb{R}^{4}$ and form the 
$4\times 4$ 
matrix $E$ whose rows are $e_{k}$, then $M=E E^{\text{tr}}$ and 
$|\Lambda|=|\det E|=\sqrt{\det{M}}$.
Using the LLL algorithm it was found that the generators of $\Sigma$ 
and $\Sigma\Sigma^{-1}$ (and hence $\Lambda$) could always be written 
as, respectively
\begin{equation}
\begin{align}
e_{k}&=2\left(a e^{i 2 k\alpha},b e^{i 2 k\beta}\right)\label{basis1}\\
e_{k}&=2\left(a(e^{i 2\alpha}-1) e^{i 2 k\alpha},b (e^{i 2\beta}-1) e^{i 
2 k\beta}\right).\label{basis2}
\end{align}
\end{equation}
To help identify the lattice geometry it was sometimes necessary to 
define a new basis $E'=S E$, where $S\in\SL{(4,\mathbb{Z})}$. The new 
Gram matrix is then given by $M'=S M S^{\text{tr}}$. 
Details of this analysis, for the seven combinations of $(n,p,q)$ in Table 
1, are provided in the appendix.

The only additional data needed in the density formula 
\eqref{rhoFormula} is the triangle area $|P|=(1/4)\sin{2\pi(p/n)}$.

The last column of Table 1 identifies which surfaces have doubly 
periodic sections.  In general, since the kernel of the homomorphism 
$\pi_{X}\colon\Lambda\to\pi_{X}(\Lambda)$ is given by the lattice 
$\Lambda_{Y}=\Lambda\cap Y$,
\begin{equation}
\Rk{\Lambda_{Y}}=\Rk{\Lambda}-\Rk{\pi_{X}(\Lambda)}.
\end{equation}
Surfaces with doubly periodic sections have $\Rk{\Lambda_{Y}}=2$, while 
$\Rk{\Lambda_{Y}}=0$
corresponds to completely quasiperiodic sections.  These are the only 
cases that occur, since (from 
\eqref{basis1} and \eqref{basis2}), $\pi_{X}(\Lambda)\simeq\mathbb{Z}[e^{2\pi i(p/n)}]$, the 
cyclotomic lattice with rank $\phi(n/\GCD{(n,p)})$.
\end{proof}

\begin{table}
\centering
{\renewcommand{\arraystretch}{2}
\begin{tabular}
{|c||c|c|c|c|c|}
\hline
$(n,p,q)$&$G/\Lambda$&$\Lambda$&$g$&$\rho$&
\begin{tabular}{c}periodic\\quasiperiodic\end{tabular}\\ \hline
$(5,1,2)$&$d_{10}$&$A_{4}$&$2$&$\frac{1}{5}\sqrt{2+\frac{2}{\sqrt{5}}}$&q\\

$(6,1,2)$&$d_{12}$&$A_{2}\times A_{2}$&$2$&$\sqrt{\frac{1}{27}}$&p\\ 

$(8,1,3)$&$d_{8}$&$D_{4}$&$2$&$\sqrt{\frac{1}{8}}$&q\\

$(10,1,3)$&$d_{10}$&$A_{4}$&$2$&$\frac{1}{5}\sqrt{2-\frac{2}{\sqrt{5}}}$&q\\ 

$(12,1,5)$&$d_{12}$&$A_{2}\times A_{2}$&$3$&$1$&q\\ 

$(12,2,3)$&$d_{24}$&$A_{2}\times Z_{2}$&$5$&$1$&p\\ 

$(12,3,4)$&$d_{24}$&$A_{2}\times Z_{2}$&$5$&$\sqrt{\frac{4}{3}}$&p\\ 

\hline
\end{tabular}
}
\caption{Properties of the seven discrete Riemann surfaces generated 
by conformal maps of right triangles}
\end{table}

\subsection{Piecewise flat surfaces and model 
sets}\label{flattenedSurfaces}

Suppose a Riemann surface $\Surf$ is deformed into another surface, 
$\Tilde{\Surf}$, not necessarily representable locally by graphs of 
holomorphic functions. As long as the deformation preserves the 
isometry group and transversality with respect to $Y$, all the 
crystallographically relevant properties of the point set $\A(x)$ 
will be maintained in the corresponding deformed point set 
$\Tilde{\A}(x)$. Kalugin's formula \eqref{KaluginFormula}, for example, makes 
this invariance explicit for the density. 

When a Riemann 
surface is generated by a triangular graph $P|Q$, deformations that preserve 
isometry and transversality are easily specified by
the map defining the fundamental graph,
$\tilde{f}:P\to Q$. We recall that in the Riemann surface, 
$\tilde{f}$ is 
holomorphic and extends to a homeomorphism on the closure $\Bar{P}$. For the 
deformed surface $\Tilde{\Surf}$ we 
continue to use the edge group $G$, defined by the geometry of the 
triangles $P$ and $Q$, to form the orbit of the fundamental graph, but 
insist only that $\tilde{f}$ is a homeomorphism. While still 
preserving isometry and transversality, we will even go a step further 
and modify the \emph{topology} of $\Tilde{\Surf}$ by defining 
$\Tilde{\Surf}$ as the orbit under $G$ of an \emph{open} graph. 
Under these circumstances, when $\Tilde{\Surf}$ is a collection of 
disconnected components, we are free to relax the condition that 
$\tilde{f}$ is a homeomorphism. In fact, we will primarily be 
interested in the case when $\tilde{f}$ is a constant map, thereby 
making each 
piece of $\Tilde{\Surf}$ flat.

Let $P|y$ represent the graph of the constant map, $\tilde{f}(P)=y$. 
We will weigh the merits of various choices of $y\in Y$ as giving optimal
``approximations'' of the map $P|Q$ defined by the conformal map 
of triangles. By choosing $y=Q_{i}$, a vertex of $Q$, we go the 
furthest in restoring partial connectedness to $\Tilde{\Surf}$. This is 
because the action of the vertex group $G_{i}$ on $P|Q_{i}$ generates 
a flat polygon (possibly stellated) composed of $|G_{i}|$ triangles. 
Thus $|G_{i}|$ surface pieces will have been ``aligned'' by this 
choice of $y$. 
With the 
exception of the $n=5$ and $n=10$ surfaces, $|G_{1}|=|G_{3}|=2 n>|G_{2}|$, 
suggesting that either one of the singular vertices is a good choice for 
$y$. 

There is another criterion, however, that applies uniformly to 
all the surfaces and even distinguishes among the two singular 
vertices.
A natural question to ask is: which value of $y\in Q$ ``occurs with 
the highest frequency'' in the graph $P|Q$? To give this question a 
proper probabilistic interpretation, we suppose that $x$ is 
sampled uniformly in $P$ and ask for the probability that 
$\|f(x)-y_{0}\|<\Delta$, where $f$ is the conformal map of triangles 
and $\Delta$ is the radius of a small disk about $y_{0}\in Q$. Since 
$f$ is conformal, this condition (for $\Delta\to 0$) is equivalent to 
$\|x-x_{0}\|<\Delta/\|f'(x_{0})\|$, where $x_{0}=f^{-1}(y_{0})$. Thus 
the probability of finding $y$ in a neighborhood of $y_{0}$ is 
maximized by minimizing $\|f'(x_{0})\|$. At a singular point 
$\|f'(x_{0})\|$ either vanishes or diverges,
and in our case vanishes for $x_{0}\to Q_{1}$ since we always have 
$p<q$. We therefore choose $y=Q_{1}$ for all of our surfaces.

The result of flattening the seven Riemann surfaces in Table 1 by this 
prescription is particularly simple.  In all cases $\Poly 
=G_{1}(P|Q_{1})$ 
is a regular $(n/p)$-gon covered $p$ times (we restore connectedness 
to $\Poly$ by including edges incident to vertex 1).
This means that if each point of $\Tilde{\A}(x)$ is counted with 
multiplicity $p$, then $\Tilde{\A}(x)$ and $\A(x)$ have the same 
density. Using Lemma \ref{groups}, we have 
\begin{equation}
\Tilde{\Surf}=\Lambda\cdot \Poly,\qquad\Lambda=\langle\Sigma\rangle,
\end{equation} 
if $r(\pi,\pi)\in R$, otherwise, 
\begin{equation}
\Tilde{\Surf}=\Lambda\:\cup\:\Lambda\, t(2 a,2 b)\, 
r(\pi,\pi)\cdot \Poly, \qquad\Lambda=\langle\Sigma\Sigma^{-1}\rangle.
\end{equation}
With the exception of the surface $(5,1,2)$, 
$r(\pi,\pi)\cdot \Poly=\Poly$, and we have the simpler description: 
\begin{equation}
\Tilde{\Surf}=\langle\Sigma\rangle\,\Poly, 
\end{equation}
since
\begin{equation}
\langle\Sigma\Sigma^{-1}\rangle\:\cup\:
\langle\Sigma\Sigma^{-1}\rangle\,t(2 a,2 b) =\langle\Sigma\rangle.
\end{equation}

Figure 1 compares the point sets $\A(x)$ and $\Tilde{\A}(x)$, given 
by the flattening process just described, of the Riemann 
surface $(5,1,2)$ of Table 1. Edges have been added to 
$\Tilde{\A}(x)$ to aid in the visualization of three tile shapes: the 
boat, star, and jester's cap. There is a one-to-one correspondence 
between the two point sets, with most pairs having quite small 
separations. In any case, we are guaranteed the separation of 
corresponding points never exceeds 1, the diameter of triangle $Q$. 
$\Tilde{\A}(x)$, the vertex set of a popular tiling model \cite{pent}, is a 
Delone set.  $\A(x)$ fails to be a Delone set because triples of 
points appear with arbitrarily short separations.  That always three 
points coalesce in this way is a signature of the order of the branch 
points of $\pi_{X}$ for this surface.  $\A(x)$, on the other hand, has 
a ``dynamical'' advantage over $\Tilde{\A}(x)$.  Seen as atoms in a 
crystal or quasicrystal, the positions $\A(x)$ evolve continuously (in 
fact analytically) with $x$ (viewed as a parameter), while the atoms 
in $\Tilde{\A}(x)$ experience discontinuous ``jumps'', and for the most 
part never move at all.  The singular loci of $x\in X$ for the two 
point sets have different dimensionalities: 0 for $\A(x)$, 1 for 
$\Tilde{\A}(x)$; the dynamics of $\A(x)$ is thus more regular also in 
this sense.  To emphasize this point, we note that if $\gamma(t)$ is 
almost any curve in $X$, then $\A(\gamma(t))$ is a regular homotopy 
(see Lemma \ref{uniformInX}).

\begin{figure}
\centerline{\epsfbox{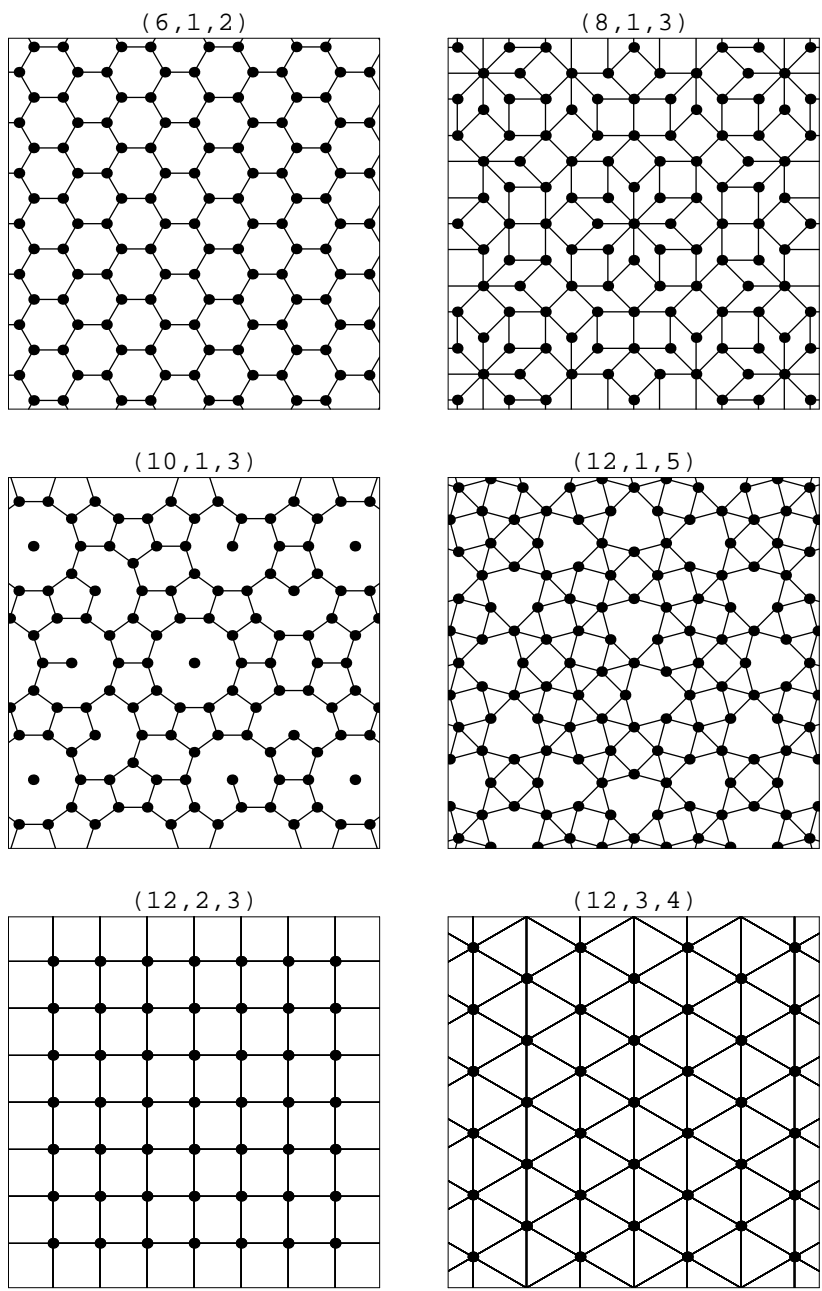}}
\caption{Point sets given by sections of flattened Riemann surfaces. 
Lines have been added to help identify tiles.}
\end{figure}

The sets $\Tilde{\A}(x)$ are examples of \emph{model sets}, as defined by 
Moody \cite{Moody}, and therefore belong to the larger family of \emph{Meyer 
sets}.  That the set $\Tilde{\A}(x)$ shown in Figure 1 can be 
organized into a finite set of tile shapes, for example, is a general 
property of model sets.  Model sets obtained by flattening the other 
six Riemann surfaces of Table 1 are shown in Figure 3.  In the three 
periodic cases, of course, even $\A(x)$ is a Meyer set and flattening 
is not necessary if that is our only goal.  We flattened these 
surfaces because the sets $\Tilde{\A}(x)$ are then particularly 
symmetric.  In the three quasiperiodic cases the sets $\Tilde{\A}(x)$ 
again organize themselves into tilings that have been discussed in the 
quasicrystal literature \cite{oct,dec,dodec}.

\section*{Acknowledgements}

I thank the Aspen Center for Physics, where a large part of this paper 
was written. Noam Elkies provided a useful suggestion for the proof of 
Lemma \ref{uniformInX}.

\newpage

\appendix
\section{Appendix}

For each of the entries in Table 1 we give below the corresponding 
Gram matrix for the generators (either \eqref{basis1} or \eqref{basis2}). From the transformed 
Gram matrices we see that the lattice geometries are in all cases 
simple root lattices ($A_{n}$, $D_{n}$, $Z_{n}$) or direct products. 
The lattices for $(5,1,2)$ and $(12,3,4)$ are obtained, 
respectively, from the lattices of $(10,1,3)$ and $(12,2,3)$ by an 
interchange of the spaces $X$ and $Y$.

\begin{align*}
(6,1,2)\qquad\qquad\qquad& &&\\
M=&
\begin{bmatrix}
6&0&-3&0\\
0&6&0&-3\\
-3&0&6&0\\
0&-3&0&6
\end{bmatrix}&&
\Lambda\simeq A_{2}\times A_{2}\qquad |\Lambda|=27\\
(8,1,3)\qquad\qquad\qquad& &&\\
M=&
\begin{bmatrix}
4&2&0&-2\\
2&4&2&0\\
0&2&4&2\\
-2&0&2&4
\end{bmatrix}
&&S=
\begin{bmatrix}
1&0&0&0\\
0&0&-1&0\\
0&-1&1&-1\\
0&0&0&-1
\end{bmatrix}\\
S M S^{\text{tr}}=&
\begin{bmatrix}
4&0&0&2\\
0&4&0&2\\
0&0&4&2\\
2&2&2&4
\end{bmatrix}&&
\Lambda\simeq D_{4}\qquad |\Lambda|=8\\
(10,1,3)\qquad\qquad\qquad& &&\\
M=\frac{1}{2}&
\begin{bmatrix}
10&5&0&0\\
5&10&5&0\\
0&5&10&5\\
0&0&5&10
\end{bmatrix}&&
\Lambda\simeq A_{4}\qquad |\Lambda|=\frac{25}{4}\sqrt{5}\\
\end{align*}
\newpage
\begin{align*}
(12,1,5)\qquad\qquad\qquad& &&\\
M=&
\begin{bmatrix}
4&3&2&0\\
3&4&3&2\\
2&3&4&3\\
0&2&3&4
\end{bmatrix}
&&S=
\begin{bmatrix}
1&-1&0&0\\
0&0&1&-1\\
1&0&-1&1\\
0&1&-1&0
\end{bmatrix}\\
S M S^{\text{tr}}=&
\begin{bmatrix}
2&1&0&0\\
1&2&0&0\\
0&0&2&1\\
0&0&1&2
\end{bmatrix}&&
\Lambda\simeq A_{2}\times A_{2}\qquad|\Lambda|=3\\
(12,2,3)\qquad\qquad\qquad& &&\\
M=\frac{1}{2}&
\begin{bmatrix}
14&3&-1&-6\\
3&14&3&-1\\
-1&3&14&3\\
-6&-1&3&14
\end{bmatrix}
&&S=
\begin{bmatrix}
1&0&1&0\\
0&1&0&1\\
1&-1&1&0\\
0&1&-1&1
\end{bmatrix}\\
S M S^{\text{tr}}=\frac{1}{2}&
\begin{bmatrix}
6&3&0&0\\
3&6&0&0\\
0&0&8&0\\
0&0&0&8
\end{bmatrix}&&
\Lambda\simeq A_{2}\times Z_{2}\qquad|\Lambda|=6\sqrt{3}
\end{align*}

\end{document}